\documentclass{amsart}

\usepackage{amssymb}
\usepackage[dvips]{epsfig}
\usepackage{graphicx}
\usepackage{color}
\usepackage[notcite,notref]{showkeys}

\newcommand{\dsp}{\displaystyle}

\newcommand{\R}{{\mathbb R}}

\newcommand{\Z}{{\mathbb Z}}
\newcommand{\T}{{\mathbb T}}

\theoremstyle{theorem}
\newtheorem{theorem}{Theorem}[section]

\newtheorem*{merci}{Acknowledgements}

\theoremstyle{remark}

\newtheorem{remark}[theorem]{Remark}

\theoremstyle{definition}

\numberwithin{equation}{section}

\begin{document}
\title[ Benjamin and related equations]{On the Benjamin and related equations}

%\title[Remarks on the KP-Benjamin and related equations]
\author[C. Klein]{Christian Klein}
\address{ Institut de Math\'ematiques de Bourgogne, UMR 5584\\   
Universit\'e de Bourgogne-Franche-Comt\'e, 9 avenue Alain Savary, 21078 Dijon
                Cedex, France}
\email{ Christian.Klein@u-bourgogne.fr}
\author{Felipe Linares}
\address{ IMPA\\ Estrada Dona Castorina 110\\ Rio de Janeiro 22460-320, RJ Brasil}
\email{ linares@impa.br}

\author{Didier Pilod}
\address{ Department of Mathematics, University of Bergen, Postbox 7800, 5020 Bergen, Norway}
\email{Didier.Pilod@uib.no}
\author{Jean-Claude Saut}
\address{Laboratoire de Math\' ematiques, Bat.307,\
Universit\' e Paris-Saclay et CNRS\\ 91405 Orsay, France}
\email{jean-claude.saut@universite-paris-saclay.fr}

\date{October18th 2023 }
\maketitle

\begin{center}
To the memory of Rafael I\'orio
\end{center}

\begin{abstract}
%\textit{Abstract}
We consider in this paper various theoretical and numerical  issues on  classical one dimensional models of internal waves with surface tension.They concern the Cauchy problem, including the long time dynamic, localized solitons or multisolitons, the soliton resolution property. We survey known results, present a few 
new ones  together with open questions and conjectures motivated by numerical simulations.

A major issue is to emphasize the differences of the qualitative behavior of solutions with those of the same equations without the capillary term.

\end{abstract}

\section{Introduction}

\subsection{The models}
This paper is concerned with one-dimensional asymptotic models for weakly nonlinear internal waves in presence of surface tension.

In  absence of surface tension, classical one-dimensional models for long, weakly nonlinear waves at the interface of a two-layer system when the lower layer is much larger than the upper one comprise the Benjamin-Ono (BO) equation derived in \cite{Ben, O}, see \cite{BLS, CGK, JCS, KS15} for more details and \cite{Pa23} for a complete rigorous justification  as an internal waves model:
\begin{equation}\label{BO}
u_t+uu_x-Hu_{xx}=0,
\end{equation}
where H is the Hilbert transform defined as  the convolution with $PV (\frac{1}{x}), $ that is  $\widehat{Hf}(\xi)=-i\text{sgn}(\xi) \hat {f}(\xi),$
(note that  here the symbol of  $H\partial_x$ is $|\xi|$ so that the full linear symbol is $-i\xi |\xi|)$, and the Intermediate Long Wave (ILW) equation derived in  \cite{KKD}
\begin{equation}\label{ILW}
u_t+uu_x+\frac{1}{\delta}u_x-\mathcal {T_\delta} u_x=0,
\end{equation}
where the \lq\lq Tilbert transform\rq\rq \, is defined  in Fourier variables by:
%$$\mathcal T_\delta u=PV\int_{-\infty}^\infty \coth\left(\frac{x-y}{\delta}\right)u(y)dy$$
%so that in Fourier variables 
$$\widehat{\mathcal{T}_\delta f}(\xi)=\xi\coth(\delta \xi)\hat{f}(\xi).$$
%The Cauchy problem for the ILW equation was shown to be globally well-posed in $H^s(\R), s>3/2$ in \cite{ABFS}.

\vspace{0.3cm}
The BO equation corresponds to the limit of an infinite lower layer, $\delta \to \infty$ in \eqref{ILW}. In this limit the symbol of $\widehat{\mathcal{T}_\delta}$ loses its regularity.

One also finds in \cite{KKD} a configuration of the layers where \eqref{ILW} should be replaced by 

\begin{equation}\label{ILWmod}
u_t+uu_x+\frac{1}{\delta_1}u_x-\mathcal T_{\delta_1}(u_x)+\frac{1}{\delta_2}u_x-\mathcal T_{\delta_2}(u_{x})=0,
\end{equation}
where $\delta_1\neq \delta_2.$

\begin{remark} Although the global well-posedness of the ILW equation in $L^2$ is now well understood (see \cite{IfSaut,CLOP}), we recall here known results on the somewhat unusual equation \eqref{ILWmod}. In   \cite{ABS}  the global well-posedness of the Cauchy problem is established in $H^s(\R), s>3/2$ by following the technique used in \cite{ABFS} for the BO and ILW equations.  The results in \cite{MV,MPV} for the usual ILW equation (global well-posedness in $H^s(\R), s\geq 1/2$ with unconditional uniqueness for $s>1/2$ and global well-posedness in $H^s(\mathbb R)$, $s>\frac14$) also apply to \eqref{ILWmod} since the symbol of this last equation satisfies Hypothesis 1 in \cite{MV,MPV}. More recently, GWP in $L^2$ both on the line and the torus were proved for \eqref{ILWmod} in \cite{CFLOP} (see Remark 1.3) by relying on the technique introduced in \cite{CLOP} for the ILW equation.

One also finds in \cite{ABS} the existence of solitary wave solutions $u(x-ct)$ for any $c>0$ by minimizing the functional on $H^{1/2}(\R),$ $J(f)=\int_{-\infty}^\infty (Lf)^2dx$ where $L$ is the Fourier multiplier defined by
$$\widehat{Lf}(\xi)=(C+m(\xi))^{1/2}\hat{f}(\xi)$$
with 
$$m(\xi)=\xi\coth(\xi \delta_1)-\frac{1}{\delta_1}+\xi\coth(\xi\delta_2)-\frac{1}{\delta_2}$$
under the constraint
$$\int_{-\infty}^\infty f(x)^3dx =1.$$

Moreover  the solitary waves are proven to decay exponentially as $|x|\to \infty.$
Note that contrary to the original  ILW equation the uniqueness (up to translation) of solitary waves to equation \eqref{ILWmod} seems to be open (see \cite{Alb1} for the ILW equation).

Concerning the stability issues, the existence of {\it orbitally stable} solitary waves is established in \cite{ABS} when $\delta_2$ is sufficiently close to $\delta_1.$ The stability in the general case is proven in \cite{Alb2}.
\end{remark} 

It is worth noticing that, similarly to the KdV equation, both \eqref{BO} and \eqref{ILW} are completely integrable (see for instance  \cite{ Fo-Ab, KSA, KAS, Ger}, \footnote{The Cauchy theory on the real line by Inverse Scattering methods is however far from being complete for those equations, in particular the soliton resolution has not been proven contrary to the KdV equation, see \cite{Schuur}. On the other hand, the periodic case was solved for the Benjamin-Ono equation, see \cite{GKT} and the references therein.  We also refer to \cite{Ger} for an {\it explicit} formula for the Benjamin-Ono equation on $\R$ and $\T$ } and the surveys \cite{JCS,KS} for this aspect). On the other hand, \eqref{ILWmod} is not known to be integrable.

We also refer to \cite{BLS} and \cite{CGK} for a rigorous justification of \eqref{BO} and \eqref{ILW} as models for long internal waves  in the sense of consistency.

\vspace{0.3cm}
All the previous equations are derived in absence of surface tension (purely gravity waves). Actually the effect of surface tension in realistic oceanic internal waves is tiny and almost negligible. Nevertheless including capillary effects in classical models of internal waves leads to interesting mathematical issues and the present paper  focusses on this aspect.

\vspace{0.5cm}
The paper is organized as follows. In the next section we present the 
derivation (often formal) of the models we will study later. Section 
3 deals with several analytical results on the Cauchy problem and 
the  solitary wave solutions, while in section 4 we present various numerical simulations that  will illustrate the previous results and help to make relevant conjectures on the long time behavior of the solutions.

We dedicate this paper to our dear colleague and friend  Rafael J. Iorio who sadly passed away on April 8, 2023.

% \subsection{Notations} \label{notations} \textcolor{blue}{is it useful?}

\section{The derivation}
We now recall the known  derivation of long, weakly nonlinear internal waves models including surface tension. 

\vspace{0.3cm}
In \cite{TBB2} Benjamin derived a version of the BO equation when the 
interface of the two layers is subject to surface tension. The modeling depends basically on the assumption that $T/g(\rho_2-\rho_1)h^2\gg 1,$ where $T$ is the interfacial surface tension, $\rho_2-\rho_1,$ the difference between the densities and $h$ the undisturbed thickness of the upper layer. This assumption is satisfied for instance when fluid densities are nearly equal.

We briefly describe here Benjamin's formal argument (see also 
\cite{Kal}). One considers  a two-fluid system made of a layer of 
fluid of thickness $h$ and density $\rho_1$ over an infinitely deep 
incompressible fluid of density $\rho_2>\rho_1,$ the less dense fluid 
being bounded from above by a rigid horizontal plane, the interface being subject of surface tension with 
coefficient $T>0.$ In this situation, the phase velocity 
$c(k)=\frac{\omega(k)}{k}$ corresponding to  plane waves $\epsilon\exp\lbrace i(\omega t-kx)\rbrace) $ reads

$$c^2(k)=\frac{g(\rho_2-\rho_1)+Tk^2}{\rho_1k\coth kh+\rho_2|k|}.$$

For waves propagating in the positive direction ($c>0$) the asymptotic approximation of $c(k)$ for very long waves ($|k|\ll1)$ gives

$$c(k)=c_0\left[1-\frac{1}{2}\frac{\rho_2}{\rho_1}h|k|+\frac{1}{2}\left\lbrace\frac{T}{g(\rho_2-\rho_1)}+\left(\frac{3}{4}\frac{\rho_2^2}{\rho^2_1}-\frac{1}{6}\right)h^2\right\rbrace k^2+O(|k|^2)\right].$$

For long waves, $|k|^2$ is much smaller that $|k|$, and thus the term in brackets can be neglected leading to the Benjamin-Ono equation. On the other hand, if it is assumed that $\rho_1$ and $\rho_2$ are very close, that is $\frac{T}{g(\rho_2-\rho_1)}h^2\gg1,$  Benjamin writes the approximate dispersion  

$$c(k)=1-\alpha |k|+\beta k^2,$$

in which  

 $$\alpha=\frac{h}{2}\rho_2/\rho_1,\quad \beta=\frac{h^2}{2}T/g(\rho_2-\rho_1).$$

This leads to the  modified BO equation, known as the Benjamin equation, again  in the notations of \cite{TBB2}\footnote{In a frame moving with velocity $c_0=\sqrt{\frac{\rho_2-\rho_1}{\rho_2}}$} :

\begin{equation}\label{BOTS}
u_t+uu_x- \alpha Hu_{xx}-\beta u_{xxx}=0,
\end{equation}

where $\beta >0$  is proportional to the surface tension coefficient. 
We have kept the parameter $\alpha>0$ to compare the effects of the 
two dispersive terms. The linear symbol is now 
$i\xi(\beta\xi^2-\alpha |\xi|)$ so that the phase velocity has  no 
longer a constant sign, contrary to the KdV or BO equations and {\it apriori }solitary waves can travel in both direction.

We refer to \cite{TBB2} and \cite{ABR} for further discussions  of the physical regime of validity of Benjamin equation. 

\begin{remark}
As aforementioned the derivation of Benjamin equation in \cite{TBB2, TBB3} is formal and a rigorous derivation in the sense of consistency is still lacking. The derivation in \cite{Kal}, is also formal but Kalish derived the following  {\it system} version of the Benjamin equation for bi-directional waves in the same physical context 

 \begin{equation}\label{BeSys}
    \left\lbrace
    \begin{array}{l}
    \dsp \zeta_t+u_x+\sigma (u\zeta)_x=0\\
    \dsp u_t+\zeta_x+\sigma uu_x-\epsilon \frac{\rho_2}{\rho_1}H\zeta_{xx}-\mu \zeta_{xxx}=0.\\
   \end{array}\right.
\end{equation}

This system appears to be a perturbation of one of the {\it abcd} Boussinesq systems derived in \cite{BCS1} corresponding to $b=d=0, a=0, c<0$ and which belongs to "Schr\"{o}dinger type" abcd systems (the dispersion is of order 2).

More precisely this (abcd) Boussinesq system is a particular case of an Euler-Korteweg system, see the survey article \cite{SX} for known results on this system. We will comment on \eqref{BeSys} in  Chapter 7.
%Actually, \eqref{BeSys} with $\epsilon=0$ is the well-posed version of the so-called Kaup-Broer-Kupperschmidt  introduced by Broer \cite{Broer} as a long wave model for surface waves and by Kaup and Kupperschmidt, \cite{Ka, Kup} in the Inverse Scattering context. This system happens to be completely integrable but this is not apparently the case of \eqref{BeSys}.

\vspace{0.3cm}
In \cite{Anh} C.T. Anh, extending the results in \cite{BLS} for internal waves without surface tension, derived rigorously in the sense of consistency from the two-layer system asymptotic systems in various regimes when surface tension is taken into account. In the Benjamin-Ono regime, he derived in the one-dimensional case the following system:

 \begin{equation}\label{ST3}
    \left\lbrace
    \begin{array}{l}
    \dsp (1+\alpha\frac{\sqrt\mu}{3}|\partial _x|)\partial_t\zeta+\frac{1}{\gamma}\partial_x((1-\epsilon\zeta) v)-\frac{\sqrt \mu}{\gamma^2}\partial_x((1-\epsilon \zeta)|\partial_x|v)\\
    \dsp+\frac{\epsilon\sqrt\mu}{\gamma^2}|\partial_x|\partial_x(\zeta v)+\left[(1-\alpha)\frac{\mu}{3\gamma}-\frac{\mu}{\gamma^3}\right]\partial_x^3v=0\\
    
    \dsp\partial_t v+(1-\gamma)\partial_x \zeta-\frac{\epsilon}{\gamma}v\partial_x v+\frac{\epsilon\sqrt\mu}{\gamma^2}\partial_x(v(|\partial_x| v))-\epsilon\sqrt \mu\nu\partial_x^3 \zeta=0,\\
       \end{array}\right.
\end{equation}
where $\epsilon\sim \sqrt \mu\ll 1$ are small parameters measuring respectively nonlinear and dispersive effects, $\gamma<1$ is the ratio of the densities of the two layers, $\nu>0$ measures the capillarity effects and $\alpha\geq 0$ is a modeling parameter. 

The classical BO system (without surface tension, see \cite{BLS}) is obtained from \eqref {ST3} when omitting  the terms of higher order $\mu, \sqrt \mu \epsilon$ and thus reads:
 \begin{equation}\label{ST4}
    \left\lbrace
    \begin{array}{l}
    \dsp (1+\alpha\frac{\sqrt\mu}{3}|\partial _x|)\partial_t\zeta+\frac{1}{\gamma}\partial_x((1-\epsilon\zeta) v)-(1-\alpha)\frac{\sqrt \mu}{\gamma^2}\partial_x|\partial_x|v =0\\
    %\dsp+\frac{\epsilon\sqrt\mu}{\gamma^2}|\partial_x|\partial_x(\zeta v)+\left[(1-\alpha)\frac{\mu}{3\gamma}-\frac{\mu}{\gamma^3}\right]\partial_x^3v=0\\
    \dsp\partial_t v+(1-\gamma)\partial_x \zeta-\frac{\epsilon}{\gamma}v\partial_x v=0.\\
       \end{array}\right.
\end{equation}
Both systems are linearly well-posed when $\alpha\geq 1.$

As noticed in \cite{BLS}, neglecting the $O(\sqrt \mu)=O(\epsilon)$ terms in \eqref{ST4}, one finds that $\zeta$ must solve a wave equation with speed $c=\sqrt{\frac{1-\gamma}{\gamma}}.$ 

Under a one-way assumption, \eqref{ST4} yields the one-parameter family of equations:

\begin{equation}\label{ST5}
\left(1+\sqrt \mu\frac{\alpha}{\gamma}|\partial_x|\right)\partial_t\zeta+c\partial_x\zeta-\epsilon\frac{3}{4}c\partial_x\zeta^2-\frac{\sqrt \mu}{2\gamma}c(1-2\alpha)|\partial_x|\partial_x\zeta=0.
\end{equation}

The usual BO equation is recovered when $\alpha=0.$

\vspace{0.3cm}
Another derivation of a system version of the Benjamin equation is made in \cite{Duran} where the following system is derived, here in its one-dimensional version
\begin{equation}\label{Duran}
    \left\lbrace
    \begin{array}{l}
    \dsp \left(1+\frac{\alpha\sqrt \mu}{\gamma}|\partial_x|\right)\partial_t\zeta+\frac{1}{\gamma}\partial _x((1+\epsilon \zeta)v)-(1-\alpha)\frac{\sqrt \mu}{\gamma^2}|\partial_x|\partial_x v=0,\\
    \dsp \partial_t v+(1-\gamma)\partial_x \zeta-\frac{\epsilon}{\gamma} v\partial_x v=T\partial_x^3v,
       \end{array}\right.
\end{equation}
where  the surface tension parameter $T$ is given by 
$$T=\frac{\sigma\mu}{g\rho_2d_1^2}.$$

In absence of surface tension, \eqref{Duran} is the Benjamin-Ono system (33) derived in \cite{BLS}. When $\mu=0,$ \eqref{Duran} is the (abcd) system in \cite{BCS1} corresponding to $a=b=d=0, c<0.$

As noticed in \cite{Duran}, \eqref{Duran} is linearly well-posed if $\alpha\geq 1.$
\end{remark}

\vspace{0.3cm}
An interesting point concerning the Benjamin equation \eqref{BOTS} is 
that it links two completely integrable equations (though being probably not integrable itself), and it thus appears to be a good candidate to test well-known features of integrable equations such as the soliton resolution.
%From now on we will take $\alpha =1.$

%\textcolor{red}{Note that in his paper \cite{Ang}, Angulo writes the equation as}

%\begin{equation}\label{Bang}
%\textcolor{red}{u_t+2uu_x-lHu_{xx}-u_{xxx}=0,\; l>0}
%\end{equation}

%\textcolor{red}{and by the change of variable $t\to -t, u\to -u$ reduces it to:}

%\begin{equation}\label{Bangbis}
%\textcolor{red}{u_t+2uu_x+lHu_{xx}+u_{xxx}=0,\; l>0}
%\end{equation}

%\textcolor{red}{He shows the existence of solitary  waves $\phi(x-ct)$ for $\gamma$ close to one, where $\gamma=\frac{l}{2\sqrt c}.$}

\vspace{0.3cm}
As Benjamin  did for the Benjamin-Ono equation, one also can take into account surface tension effects in the ILW equations \eqref{ILW}, \eqref{ILWmod}, yielding respectively

\begin{equation}\label{ILWTS}
u_t+uu_x+\frac{1}{\delta}u_x-\mathcal T_\delta(u_{x})-\beta u_{xxx}=0,
\end{equation}

\begin{equation}\label{ILWmodTS}
u_t+uu_x+\frac{1}{\delta_1}u_x-\mathcal T_{\delta_1}(u_{x})+\frac{1}{\delta_2}u_x-\mathcal T_{\delta_2}(u_{x})-\beta u_{xxx}=0,
\end{equation}

From now on we will refer to  \eqref{ILWmodTS} as the {\it ILW-Benjamin (ILW-B) equation.} As does the Benjamin equation this equation connects two integrable equations, the KdV and the ILW equations yielding interesting issues on the long time behavior of solutions.

\begin{remark}
We will not comment more on the  derivation of  equations \eqref{ILWTS}, \eqref{ILWmodTS} but the formal approach of Benjamin provides  a justification. We refer again to \cite{Anh} for a rigorous derivation (in terms of consistency) of system versions (in both one and two spacial dimensions) of the BO and ILW equations following the approach of \cite{BLS} for the purely gravity waves. 

As aforementioned their one-way propagation versions do not reduce to  the Benjamin or ILW-B equations, they contain higher order terms that are neglected in the Benjamin or ILW-B equations. The situation here is somewhat similar to the formal derivation of the Kawahara equation which modifies the KdV equation 
\begin{equation}\label{KdVKawa}
\eta_t+c_0\eta_x+\frac{3c_0}{2h}\eta\eta_x+\frac{1}{2}c_0h^2(\frac{1}{3}-\tau)\eta_{xxx}=0,
\end{equation}
when the Bond number $\tau$ is very close to $\frac{1}{3}.$ This leads (see for instance  \cite{HS}) to the Kawahara equation
\begin{equation}\label{Kawa}
\eta_t+c_0\eta_x+\frac{3c_0}{2h}\eta\eta_x+\frac{1}{2}c_0h^2(\frac{1}{3}-\tau)\eta_{xxx}+\frac{1}{4}c_0h^4(\frac{14}{90}-\tau^2)\partial_x^5\eta=0,
\end{equation}
where higher order nonlinear terms are neglected.
\end{remark}

\vspace{0.5cm}
The equations \eqref{ILWTS}, \eqref{BOTS} can be seen as perturbations of the Korteweg-de Vries equation (KdV)
\begin{equation}\label{KdV}
u_t+uu_x-\beta u_{xxx}=0.
\end{equation}
They can alternatively be written in the form
\begin{equation}\label{eqL}
u_t +uu_x-\mathcal Lu_x-\beta u_{xxx}=0,
\end{equation} 
where $\mathcal L$ is defined in Fourier space by
$$\widehat{\mathcal L f}(\xi)= p(\xi) \widehat{f}(\xi),$$
with $p(\xi)=|\xi|$ for equation \eqref{BO},  $p(\xi)=p_\delta(\xi)=\xi\coth(\delta \xi)-\frac{1}{\delta}$ for equation \eqref{ILW}  and
$p_{\delta_1,\delta_2}(\xi)= \xi\coth(\delta_1 \xi)+\xi\coth(\delta_2 \xi)-\frac{1}{\delta_1}-\frac{1}{\delta_2}$  for equation \eqref{ILWTS}.

To make more clear the interaction of the two dispersive  terms we introduce a parameter $\alpha>0$ and rewrite the equations:
\begin{equation}\label{BOBmod}
u_t+uu_x-\alpha Hu_{xx}-\beta u_{xxx}=0,
\end{equation}
\begin{equation}\label{ILWBmod}
u_t+uu_x+\alpha\left(\frac{1}{\delta}u_x-\mathcal T_\delta(u_{x})\right)-\beta u_{xxx}=0.
\end{equation}
and 
\begin{equation}\label{ILWmodBmod}
u_t+uu_x+\alpha\left(\frac{1}{\delta_1}u_x-\mathcal T_{\delta_1}(u_{x})+\frac{1}{\delta_2}u_x-\mathcal T_{\delta_2}(u_{x})\right)-\beta u_{xxx}=0,
\end{equation}

Note, however, that contrary to the KdV equation the phase velocity $c(\xi)=\frac{\omega(\xi)}{\xi}$ where $\omega$ is the dispersion relation has a non constant sign. For instance, for equation \eqref{BOBmod} one finds  $c(\xi)=\beta\xi^2-\alpha |\xi|$  and $c(\xi)=\beta \xi^2-\alpha(\xi\coth(\delta \xi)-\frac{1}{\delta})$ for \eqref{ILWBmod}. 

This fact has an influence on the admissible velocities of solitary wave solutions, see \cite{TBB2, TBB3} and below.

Another important point is that none of those equations appear to be completely integrable contrary to  the KdV, Benjamin-Ono or Intermediate Long Wave equations. They are thus 
very good examples to test the expected soliton resolution conjecture.

% Actually, one easily checks via Pohojaev type arguments and the use of the interpolation inequality 

%$$|D^{1/2}u|^2_2\leq \frac{\epsilon}{2}|u|^2_2+\frac{1}{2\epsilon}|u_x|^2_2,\; \epsilon >0,$$

%that no non trivial solitary waves of \eqref{BOTS} with positive velocity $c$ exist  when $c<\frac{1}{5\beta}.$

%A similar result holds for equations \eqref{ILWTS} and \eqref{ILWmodTS}.

%Concerning the Cauchy problem, Linares has proven in \cite{Lin2} that \eqref{BOTS} is globally well-posed for initial data in $L^2(\R)$ and also, in a periodic setting, in $L^2(\T).$  It is very likely that such a result extends to \eqref{ILWTS} and \eqref{ILWmodTS}. \textcolor{red}{should we add something here?}.

%Well-posedness below $L^2(\R)$ was obtained in \cite{CGX, KOT, LW} leading to global well-posedness in $H^s(\R), s\geq -3/4.$. 

%We refer to \cite{ABS, ABR, Ang} and to Subsection 1.2 below for existence and stability properties of solitary waves to \eqref{BOTS}, \eqref{ILWTS}, \eqref{ILWmodTS}.

% Note that due to the lack of regularity of the symbol of \eqref{BOTS} the associated (non explicit) solitary waves have only an algebraic decay while the decay is exponential for the two other equations.

%While 
%the global well-posedness is well established for those equations, the global behavior of solutions is unknown and one scope of this paper is to give in particular  numerical evidence of the soliton resolution.

\begin{remark}
There are interesting Kadomtsev-Petviashvili  (KP) versions of those one-dimensional models, that is taking into account weak transverse effects, which should be useful to study the transverse stability of the one-dimensional solitary waves.

Proceeding as in the formal derivation of the classical Kadomtsev-Petviashvili equation in \cite{KaPe} (see for instance \cite{Kim3, KA} for a formal derivation of \eqref{BOBKP} from the two-layer system for the BO case) one obtains respectively 
\begin{equation}\label{BOBKP}
u_t+uu_x-\alpha  Hu_{xx}-\beta u_{xxx}+\partial_x^{-1}u_{yy}=0,
\end{equation}
\begin{equation}\label{ILWBKP}
u_t+uu_x+\alpha\left(\frac{1}{\delta}u_x-\mathcal 
T_\delta(u_{x})\right)-\beta u_{xxx}+\partial_x^{-1}u_{yy}=0,
\end{equation}
and 
\begin{equation}\label{ILWmodBKP}
u_t+uu_x+\alpha\left(\frac{1}{\delta_1}u_x-\mathcal 
T_{\delta_1}(u_{x})+\frac{1}{\delta_2}u_x-\mathcal T_{\delta_2}(u_{x})\right)-\beta u_{xxx}+\partial_x^{-1}u_{yy}=0.
\end{equation}

We refer to a subsequent work for a study of these equations that are mathematically interesting since they have both a focusing and defocusing aspect (depending on frequencies). The KP versions of the Benjamin-Ono and ILW equations are studied in \cite{LPS2} but many issues remain open.

%Proceeding as in \cite{Lannes1} (see also \cite{La, LS}) for the classical KP-II equation the derivation can probably be justified (in the sense of consistency) starting from the fully two-dimensional system in \cite{BLS} taking into account a surface tension extra term.

%\begin{remark}
%Contrary to the KP versions of the classical BO or ILW equations (without surface tension), which are of KP-II type, {\it ie} defocusing (see \cite{LPS2}) one has here KP-I type equations, which are in some sense focusing at least for large frequencies. Note however that for small frequencies \eqref{BOBKP}, \eqref{ILWBKP}, \eqref{ILWmodBKP} are reminiscent of the KP-II BO or KP-II ILW equations that were formally derived in \cite{AS, CC}  and this should have an effect on the global behavior of solutions.
%\end{remark}

Note also that due to the second order nonlocal dispersive operator, the KP-Benjamin type  equations do not possess the scaling invariance property 
$$u_\lambda(t,x,y)=\lambda^2u(\lambda^3t,\lambda x,\lambda^2y),\; \lambda>0$$
of the usual KP equation.
\end{remark}

\vspace{0.5cm}
%The paper is organized as follows. In a  first Chapter we continue  to 
%survey  the known results (and their  relatively easy extensions to the ILW-Benjamin equation ) on the Benjamin and related equations. The next Chapter present new results on the Benjamin and ILW-Benjamin equations . In Chapter 4 we give some results on the system version \eqref{BeSys} of the Benjamin equation.

%In a final Chapter, various numerical simulations will illustrate the previous results and help to make relevant conjectures on the long time behavior of the solutions.

%\subsection{Notations} \label{notations}

\section{Analysis}
The aim of this section is to review the known results on the Cauchy problem and on solitary wave solutions  for the Benjamin equation and to extend some of them   to the ILW-Benjamin equation.

\subsection{The Cauchy problem}

\vspace{0.3cm}
Concerning the Cauchy problem for the Benjamin equation, Linares has proven in \cite{Lin2} that \eqref{BOTS} is globally well-posed for initial data in $L^2(\R)$ and also, in a periodic setting, in $L^2(\T).$  
Well-posedness below $L^2(\R)$ was obtained in \cite{CGX, KOT, LW} leading to global well-posedness in $H^s(\R), s\geq -3/4.$ 

We now briefly explain how these results extend to \eqref{ILWTS}, \eqref{ILWmodTS}. By writing these equations on the form (1.3) in \cite{MV}, we observe that the symbol  of the respective linear symbols of these equations 
 $p_{ILW-B}=\beta \xi^3-(\coth{\delta xi}\xi-\frac1{\delta})\xi$ and $p_{mILW-B}=\beta \xi^3-(\coth{\delta_1xi}\xi-\frac1{\delta_1})\xi-(\coth{\delta_2xi}\xi-\frac1{\delta_2})\xi$ satisfy Hypothesis 1 in \cite{MV} (this is a direct consequence of Lemma 2.1 in \cite{MV}). Therefore, we deduce from Theorem 1.5 and Corollary 1.10 in  \cite{MV} that the Cauchy problems associated with \eqref{ILWTS} and \eqref{ILWmodTS} are globally well-posed in $L^2(M)$ with $M=\mathbb R$ or $\mathbb T$.

%We refer to \cite{ABS, ABR, Ang} and to Subsection 1.2 below for existence and stability properties of solitary waves to \eqref{BOTS}, \eqref{ILWTS}, \eqref{ILWmodTS}.

% Note that due to the lack of regularity of the symbol of \eqref{BOTS} the associated (non explicit) solitary waves have only an algebraic decay while the decay is exponential for the two other equations.

\vspace{0.3cm}
While 
the global well-posedness is well established for those equations, the global behavior of solutions is unknown and challenging. Being in some sense an \lq\lq interpolate" equation between two integrable ones, one might expect that  the soliton resolution holds. A first preliminary step towards this conjecture would be  the existence and stability of multi-soliton solutions but this seems to be    an open question as we will discuss later on.
 %For instance, unique continuation properties were established in \cite{Cun, Ur}.

\begin{remark}
Many papers have been devoted to controllability properties of the KdV or Benjamin-Ono equations, see for instance \cite{LLR, Li-Ro} for the BO case. Concerning the Benjamin equation, the controllability and stabilization were proven in \cite{PV2} for the linearized  equation and in \cite{PV} for the  equation itself, both on a periodic domain.

%We briefly indicate here how some of those  results can be extended to \eqref{ILWTS} and  \eqref{ILWmodTS}. 

%\textcolor{red}{Felipe, do you think that the extension of the $L^2$ theory is straightforward?}

%\vspace{0.3cm}
%Concerning the Cauchy problem for the KP version, the easy results known for the KP-I equation, see for instance \cite{KS}, that is local well-posedness in a suitable subspace of $H^s(\R^2), s\geq 3$ or global existence of weak solutions in the energy space $Y$ extend easily to equations \eqref{BOBKPbis}, \eqref{ILWBKPbis},\eqref{ILWmodBKPbis}.

%The deeper results concern mainly \eqref{BOBKP} which has been first   studied mathematically  in \cite{Z1, Z2}. Zaiter established in \cite{Z2} that \eqref{BOBKP} is quasilinear in the sense of \cite{MST3}, {\it ie}  that the Cauchy problem cannot be solved by an iterative method applied to the integral Duhamel formulation, for initial data in Sobolev spaces. She also extended to \eqref{BOBKP} the results of \cite{MST1} concerning the zero mass constraint.

%Those results extend easily to equations \eqref{ILWBKPbis}, \eqref{ILWmodBKPbis}.

%Very recently, it was established in \cite{CWL} that the Cauchy problem for  \eqref{BOBKP} is globally well-posed in the energy space

%$$Y=\lbrace u\in L^2(\R^2); \partial_x u, \partial_x^{-1}\partial_yu \in L^2(\R^2)\rbrace,$$

%extending thus a similar result established in \cite{IKT} for the KP-I equation.

%The global well-posedeness in the energy space for \eqref{BOBKPbis} was established in \cite{CWL}. 

It would be interesting to extend those results to \eqref{ILWTS} and  \eqref{ILWmodTS}. 
\end{remark}

%\textcolor{red}{TO BE DONE}

\subsection{Solitary waves}

We refer to \cite{ABS, ABR, Ang}  for existence and stability properties of solitary waves to \eqref{BOTS}, and to the end of the present chapter for extensions to \eqref{ILWTS}, \eqref{ILWmodTS}.

 Note that due to the lack of regularity of the symbol of \eqref{BOTS} the associated (non explicit) solitary waves have only an algebraic decay while the decay is exponential for the two other equations.

We first review and complete   known results on  the existence and properties of solitary waves, starting by the BO case that is for the equation \eqref{BOBmod}. A solitary wave solution $u(x,t)=v(x-ct)$  of \eqref{BOBmod} satisfies the equation
\begin{equation}\label{SWBOBmod}
-cv-\beta v_{xx}-\alpha|D|v+\frac{v^2}{2}=0
\end{equation}
yielding the energy equality
\begin{equation}\label{E1}
\int_{\R}\left(-cv^2+\beta v_x^2-\alpha |D^{1/2}v|^2+\frac{v^3}{2}\right)dx=0.
\end{equation}
 A Pohojaev identity is obtained by multiplying \eqref{SWBOBmod} by $xv_x$ and integration (justification by the usual smoothing and truncation process). We observe by Plancherel and integration by parts that
$$\int_\R xv_x|D|vdx=0$$
and obtain
\begin{equation}\label{E2}
\int_{\R}\left(\frac{c}{2}v^2+\frac{\beta}{2}v_x^2-\frac{v^3}{6}\right)dx=0
\end{equation}
and by combination of those identities
\begin{equation}\label{E3}
\int_{\R}\left(\frac{c}{2}v^2+\frac{5\beta}{2}v_x^2-\alpha|D^{1/2}v|^2\right)d x=0,
\end{equation}
that is in Fourier variables
\begin{equation}\label{E4}
\int_{\R}\left(\frac{c}{2}+\frac{5\beta}{2}|\xi|^2-\alpha|\xi|\right)|\hat{v}|^2d\xi=0.
\end{equation}

 We recover here the well-known fact that for the pure KdV equation ($\alpha=0$) there exists a unique depression solitary wave with velocity c  strictly negative (mind the minus sign in front of $u_{xxx}!$)  while for  the pure Benjamin-Ono equation ($\beta=0$) there exists a unique elevation solitary wave with  a strictly positive velocity.
 
 In the present case, \eqref{E4} implies that no solitary wave with 
 large velocities exists, that is  when $c>\frac{\alpha^2}{5\beta}$ and a natural conjecture is that solitary waves solutions of \eqref{SWBOBmod}
exist with velocities $c\in ( -\infty, \frac{\alpha^2}{5\beta})$ .
%\textcolor{red}{USELESS?}

\vspace{0.4cm}
The existence of solitary wave solutions of the Benjamin equation was sketched in \cite{TBB2} by minimizing the Hamiltonian with fixed $L^2$ norm.
Using a Leray-Schauder degree theory, Benjamin \cite{TBB3} proved the existence of solitary waves (and also of space periodic travelling waves)  
for small values of the parameter $\alpha/\sqrt{\beta c}$, by perturbation of a KdV soliton.

Chen and Bona, \cite{Ch-Bo}, using a concentration-compactness approach, proved the existence of solitary waves of velocity c provided $\min_{x\in \R}\lbrace x^2-\alpha |x|+c\rbrace >0,$ which implies in our notations $c<-\frac{\alpha^2}{4\beta}.$

 The later paper contains also a rigorous asymptotics of the solitary waves following the general theory of Bona and Li, \cite{Bo-Li}, namely as the (explicit) soliton of the Benjamin-Ono equation, any solitary wave solution of the Benjamin equation has the asymptotics 
 
 $$\lim_{x\to \pm \infty} x^2\phi(x)=C,$$
 for some constant $C\in \R, C\neq 0.$

 In \cite{Ang} Angulo minimizes the energy functional
 $$E(\psi)=\frac{1}{2}\int_\R\left(\beta \psi_x^2-\psi H\psi_x-\frac{1}{3}\psi^3\right)dx$$
 under the constraint 
 $$\int_{-\infty}^\infty \psi^2 dx=\lambda>0.$$
 This yields, for any $\lambda>0$ small enough, the  existence of  a nonempty set $G_\lambda$ of solitary waves with positive stability which is orbitally stable.
 \begin{remark}
 In \cite{Ang2}, Angulo proved the instability of solitary waves of the {\it generalized} Benjamin equation, extending similar known results for the generalized KdV equation.
 \end{remark}
 
 To summarize the above results, existence and stability of solitary waves for the Benjamin equation is known in the velocity range
 $$-\infty<c<-\frac{\alpha^2}{4\beta}$$
 and non existence in the range 
 $$c>\frac{\alpha^2}{5\beta}.$$
 Existence for velocities satisfying 
 $$-\frac{\alpha^2}{4\beta}<c<\frac{\alpha^2}{5\beta}$$
 is therefore an open problem that will be addressed numerically  in the present paper.  Note again that the possibility of existence of solitary waves with positive velocities is due to the BO term in the Benjamin equation. 
 %\textcolor{red}{DIDIER OR NUMERICAL SIMULATIONS?}
 
 \begin{remark}
 A stationary solution of the Benjamin equation should satisfy
 \begin{equation}\label{Ben-stat}
 -\beta v_{xx}-\alpha|D|v+\frac{v^2}{2}=0.
 \end{equation}
 While no such non trivial solutions exist by \eqref{E4} in the pure KdV or Benjamin-Ono case, the Pohajaev type argument does not exclude the possibility of existence of steady solutions for the Benjamin equation. This issue appears to be open.
 \end{remark}
 
 \begin{remark}
 Concerning the qualitative properties of the Benjamin solitary waves we already mentioned that Bona and Chen \cite{Ch-Bo} gave the algebraic decay rate of the solitary waves. They also conjecture that the solitary waves are even, have a finite number of oscillations and decay in a monotone way at $\pm \infty.$  This issue appears to be open.
 \end{remark}
 
 \begin{remark}
 An interesting open issue is that of the uniqueness (up to the trivial symmetries) of the solitary waves. Actually the techniques used in \cite{AmTo, AmTo2, Alb1} does not work because of the non-positivity of the symbol $-\alpha|\xi|+\beta \xi^2.$
 
 Also, since contrary to the KdV or Benjamin-Ono equations the Benjamin equation is not integrable, the possible soliton decomposition property is an important open problem.
 \end{remark}

 \vspace{0.3cm}
 We refer to \cite{DDM,DDM2} and to section 4 for numerical simulations of solitary waves of the Benjamin and Benjamin type equations.
 
% \begin{remark}
%The uniqueness of traveling wave solutions to the Benjamin or ILW-Benjamin equations seems to be an open problem contrary to the case of the Benjamin-Ono or ILW equation (see \cite{AmTo, AmTo2, Alb1}).
%\textcolor{red}{This is an interesting problem; I asked John Albert about it. He told me that the method he used to prove uniqueness of the BO or ILW solitary wave does not work here by lack of a positivity property.}
%\end{remark}

\begin{remark}

By an implicit functions argument from the cnoidal wave solutions of the KdV equation the existence of an analytic curve of periodic solitary waves of the Benjamin equation  was proven in \cite{AA} 
together with their orbital stability. We briefly sketch  here how to extend those results to the ILW-B equation written as 
%\textcolor{red}{JEAN-CLAUDE OR IN A SUBSEQUENT WORK?}

\begin{equation}\label{ILW-Bnew}
u_t+uu_x+u_{xxx}-l\mathcal T_\delta(u_{xx})=0
\end{equation}
\end{remark}
where $l\in \R.$

When acting on $2L-$ periodic functions of $x$,   the nonlocal operator $ \mathcal T_\delta$ is defined by (see \cite{AFSS, ABFS})
$$\mathcal T_\delta f(x)=\frac{1}{\delta} f_x+i\sum_{k\in \Z^*}\coth\left(\frac{k\pi\delta}{L}\right)f_ke^{ik\pi x/L}.$$
We are looking for real even periodic traveling wave solutions of the form $u(x,t)=\phi(x-ct)$ of \eqref{ILW-Bnew}, where $\phi=\phi_{l,c}$ satisfies the equation:
\begin{equation}\label{cnoid}
\phi''-l\mathcal  T_\delta \phi+\frac{1}{2}\phi^2-c\phi=A_{\phi_{l,c}},
\end{equation}
where $A_{\phi_{l,c}}$ is  a constant of integration which will be different from zero. Actually if $ \phi_{l,c}$ has mean zero over $[-T,T],$ the constant of integration satisfies
$$A_{\phi_{l,c}}=\frac{1}{4L}\int_{-L}^L\phi^2_{l,c}(x)dx.$$
When $l=0,$ $\phi_{0,c}$ represents the well known periodic solution of the equation 
$$\psi''_c+\frac{1}{2}\psi_c^2-c\psi_c=A_{\psi_c},$$
the so-called cnoidal wave of the KdV equation.

It is very likely that the implicit function approach used in \cite{AA} to prove existence of a continuous curve of periodic solitary waves of the Benjamin equation can be used to get a similar result for the ILW-Benjamin equation and thus solutions to equation \eqref{cnoid}.

\vspace{0.3cm}
We now comment on  the existence of solitary waves to the Benjamin-ILW equation. A natural approach is to look for solutions to the following variational problem:
%\textcolor{red}{JEAN-CLAUDE OR SUBSEQUENT WORK?}
\begin{equation}\label{var}
 \inf\lbrace  \mathcal{H}(u); \int_\R u^2dx =constant\rbrace,
\end{equation}
where 
$$\mathcal {H}(u)=\frac{1}{2}\int_{\R^2}\lbrack \beta u_x^2-|T_\delta^{1/2}u|^2+\frac{1}{3}u^3\rbrack,$$
as sketched by Benjamin \cite{TBB2} and implemented in \cite{Ang} for the Benjamin equation. 

Another approach used in \cite{Ch-Bo} for the Benjamin equation would consist in considering the minimization problem for 
$\lambda>0:$
$$\inf\{ J(u), u\in H^1(\R) ;\ \int_\R u^3(x)dx =\lambda\},$$
where
$$J(u)=\int_\R u(c+\mathcal L)udx,$$
$\mathcal L$ being the operator with symbol $|\xi|^2-l(\xi\coth(\xi \delta)-1/\delta).$
We refer to a subsequent work for an implementation of those methods.

%The known results concern mainly \eqref{BOBKP} which has been first   studied mathematically  in \cite{Z1, Z2}. Zaiter established in \cite{Z2} that \eqref{BOBKP} is quasilinear in the sense of \cite{MST3}, {\it ie}  that the Cauchy problem cannot be solved by an iterative method applied to the integral Duhamel formulation, for initial data in Sobolev spaces. She also extended to \eqref{BOBKP} the results of \cite{MST1} concerning the zero mass constraint.

%Very recently, it was established in \cite{CWL} that the Cauchy problem for  \eqref{BOBKP} is globally well-posed in the energy space

%$$Y=\lbrace u\in L^2(\R^2); \partial_x u, \partial_x^{-1}\partial_yu \in L^2(\R^2)\rbrace,$$

%extending thus a similar result established in \cite{IKT} for the KP-I equation.
\vspace{0.3cm}
The symbol $\xi\coth(\xi \delta)-1/\delta$ of the ILW-Benjamin equation being smooth, contrary to that of the Benjamin equation, one can expect {\it exponential} decay of  the solitary waves of the ILW-Benjamin equation.
% and actually :

%\begin{proposition}
%Any solitary wave $\phi$ of the ILW-Benjamin equation is smooth, $\phi\in C^\infty(\R),$  and satisfies the decay property:
%\end{proposition}

%\begin{proof}
A natural idea is to follow the strategy of Bona-Li, \cite{Bo-Li}, writing the equation of a solitary wave $\phi$ as 

$$-\phi''+\mathcal{T}_\delta \phi+c\phi=\frac{1}{2}\phi^2$$

or 

$$\phi=\frac{1}{2}K\star \phi^2,$$

where $K$ is the inverse Fourier transform of

Concerning the decay property, the first step would be  to obtain the asymptotic behavior of the Fourier transform of the function 

$$p(\xi)=\left(c+\beta\xi^2-\left(\xi \coth(\delta \xi)-\frac{1}{\delta}\right)\right)^{-1}.$$

%\textcolor{red}{SUBSEQUENT WORK?}
%\end{proof}

\vspace{0.3cm}

\subsubsection{Symmetry and uniqueness of solitary waves} 
The profile  of Benjamin solitary waves appears to be symmetrical, see the simulations below. 

% The figures represent a perturbation of the solitary wave $Q$ by 
% 
% $0.05\exp(-(x-x_0)^2),$ with $x_0=\pm 1,0.$

% \begin{figure}[!b]
%     \begin{centering}
%     \includegraphics[width=12cm]{./SW-Ben1}
%     \end{centering}
% \end{figure}
% 
% \begin{figure}[!b]
%     \begin{centering}
%     \includegraphics[width=12cm]{./SW-Ben2}
%     \end{centering}
% \end{figure}
% 
% \begin{figure}[!b]
%     \begin{centering}
%     \includegraphics[width=12cm]{./SW-Ben3}
%     \end{centering}
% \end{figure}

While the solitary waves of the KdV and Benjamin-Ono equation are known to be unique (see \cite{AmTo} for the delicate case of the BO equation) no such result seems to be known for the Benjamin equation.

\begin{remark}

1. The existence and possible stability of multi-solitons is an interesting open question for the Benjamin type equations.

%\subsection{Soliton resolution}

%The Benjamin equation being not integrable, there is no theoretical results on the soliton resolution. However this phenomena is confirmed by the following numerical simulations.

%\subsection{Scattering of small solutions} 

2. The non existence of small solitary waves suggests the scattering of small solution which to our knowledge is an open question.
\end{remark}
.

%Similarly,

%$$G_2(x,t)=e^{\frac{i}{3\beta}(\frac{2t}{9\beta}-x)}\int_0^\infty e^{i\eta(x-\frac{t}{3\beta})}e^{it\eta^3} d\eta =t^{-1/3}e^{\frac{i}{3\beta}(\frac{2t}{9\beta}-x)}\int_{-\infty}^0 e^{it^{-1/3}\eta(x-\frac{t}{3\beta})}e^{i\eta^3}d\eta$$

%We are thus reduced to compute the asymptotics of the functions:

%$$A_1(x)=\int_{-\infty}^0 e^{ix\xi}e^{i\xi^3}d\xi,$$

%and
%$$A_2(x)=\int_0^\infty e^{ix\xi}e^{i\xi^3}d\xi.$$

\subsection{Unique continuation}

In this subsection we are concerned with unique continuation principles (UCP) satisfied for solutions to the equations considered in this paper. 

In \cite{KPV-ucp} Kenig, Ponce and Vega proved the following linear result for generalizations of the BO equation. 
Suppose that $k, j\in {\Z}^{+}\cup \{0\}$ and that
\begin{equation*}
a_m:\R\times [0,T] \mapsto \R, \hskip5pt m=0,1,\dots,k,  \hskip5pt\text{and}\hskip5pt b:\R\times [0,T] \mapsto \R
\end{equation*}
are continuous functions with $b(\cdot)$ never vanishing on $(x,t)\in\R\times[0,T]$, and consider the IVP
\begin{equation}\label{w-linear}
\begin{cases}
\partial_t w-b(x,t)H\partial_x^jw+\underset{m=0}{\overset{k}{\sum}} a_m(x,t)\, \partial_x^mw=0,\\
w(x,0)=w_0(x).
\end{cases}
\end{equation}

\begin{theorem}[Theorem 1.4 \cite{KPV-ucp}]\label{ucp-linear}
Let  $w\in C([0,T]:H^s(\R))\cap C^1((0,T):H^{s-2}(\R))$, $s>\max\{k\,;j\}+1/2$ be a real solution to the IVP \eqref{w-linear}. If there exists an
open set $\Omega\subset \R\times[0,T]$ such that
\begin{equation}\label{w-linear-1}
w(x,t)=0,\hskip10pt (x,t)\in\Omega,
\end{equation}
then 
\begin{equation}\label{w-linear-2}
w(x,t)=0,\hskip10pt (x,t)\in \R\times[0,T].
\end{equation}
\end{theorem}

As pointed out in \cite{KPV-ucp} the following result can be seen as a corollary of Theorem \ref{ucp-linear}. 
\begin{theorem}[Theorem 1.1 \cite{KPV-ucp}]\label{ucp-bo}
Let  $u_1, u_2$ be real solutions to the IVP associated to the BO equation for  $(x,t)\in \R\times[0,T]$ such that
$u_1, u_2\in  C([0,T]:H^s(\R))\cap C^1((0,T):H^{s-2}(\R))$, $s>5/2$. If there exists an
open set $\Omega\subset \R\times[0,T]$ such that
\begin{equation}\label{bo-1}
u_1(x,t)=u_2(x,t),\hskip10pt (x,t)\in\Omega,
\end{equation}
then 
\begin{equation}\label{bo-2}
u_1(x,t)=u_2(x,t),\hskip10pt (x,t)\in \R\times[0,T].
\end{equation}
In particular, if $u_1$ vanishes in $\Omega$, then $u_1\equiv 0$.
\end{theorem}

An application of Theorem \ref{ucp-linear} to the difference of two real solutions $u_1, u_2$ of the Benjamin equation \eqref{BOBmod}  implies
that the result in Theorem \ref{ucp-bo}  with $s>7/2$, holds for the IVP associated to the Benjamin equation \eqref{BOBmod}.

In a similar direction, solutions associated to IVP of the ILW 
equation were proved to have the same properties as the solutions of the BO equation
mentioned above. For details see Theorem 1.8  and the remark 1.9 in \cite{KPV-ucp}. We can deduce then that solutions of the IVP associated to
the equations  \eqref{ILWBmod} and \eqref{ILWmodBmod} for $s>7/2$ 
satisfy a Unique Continuation Principle (UCP) as the one described in Theorem \ref{ucp-bo}.

A stronger UCP has been obtained for the Benjamin equation \eqref{BOBmod}.  Consider the functional space 
$Z_{s,r}(\R):= H^s(\R)\cap L^2(|x|^{2r}\,dx)$. We recall that for solutions of the BO equation Iorio  (\cite{I}) proved that
if $u\in C([0,T], Z_{2,2}(\R))$, and there exist
three distinct times $t_1, t_2, t_3$ such that the solutions satisfy
\begin{equation}
u(\cdot, t_j) \in Z_{4,4}(\R), \hskip10pt j=1,2,3,
\end{equation}
then $u\equiv 0$. This result is sharp as it was shown in \cite{FLP}. A similar result was established in
\cite{Ur}  for solutions of the Benjamin equation in weighted spaces. 
\begin{theorem}[Theorem 3 in \cite{Ur}] Suppose that $u\in C([0,T]: \dot{Z}_{7,7/2^{-}})$ is a solution of the IVP associated to the Benjamin equation.
If there exist three different times $t_1, t_2, t_3\in [0,T]$ such that $u(\cdot, t_j) \in \dot{Z}_{7,7/2^{-}}$ for $j=1, 2, 3$, then $u(x,t)\equiv 0$. Here
$\dot{Z}_{s,r}:= \{ f\in Z_{s,r} : \hat{f}(0)=0\} $.
\end{theorem}
As far as we know no such result is available for solutions of the ILW-Benjamin equations \eqref{ILWBmod} and \eqref{ILWmodBmod}.

%\begin{remark}

%\textcolor{red}{SUBSEQUENT WORK TOGETHER WITH THE KP CASE?}
%We comment  here  on  the  system version  \eqref{BeSys} of the Benjamin equation .  It has the Hamiltonian structure

%$$\partial_tU+\mathcal J\mathcal H (U)=0$$

%where 

%$$U=\begin {pmatrix} \zeta\\ u\end{pmatrix},\quad \mathcal J=\begin{pmatrix}0&\partial_x\\\partial _x&0\end{pmatrix},\quad\mathcal H=\left(u\zeta+\frac{\sigma}{2}u^2\zeta+\frac{\mu}{2}\zeta_x^2-\frac{\epsilon\rho_2}{2\rho_1}|D^{1/2}\zeta|^2\right) $$

%Unfortunately, as for the associated (abcd) system (corresponding to $\epsilon =0$) , the conservation of the Hamiltonian does not provide the apriori estimates necessary to obtain the global existence of small solutions.

%As aforementioned, \eqref{BeSys} is a perturbation of the Boussinesq "Schr\"{o}dinger type" (abcd) system in the terminology of \cite{BCS1} corresponding to $a=b=d=0, c<0$ for which long time existence (that is on time scales of order $O(1/\epsilon)$) of the Cauchy problem has been obtained in \cite{SWX}. 

%We briefly indicate how to extend this result to \eqref{BeSys}. 
%An other issue is that long time existence, that is on time scales of order $O(1/\epsilon).$ This was established in \cite{SWX} for the  corresponding Boussinesq system. 
%We briefly show how to extend this result to the case of system \cite{BeSys}.
%
\section{Numerical simulations}
We illustrate here  qualitatively the behavior of solutions to the 
Benjamin equation aiming to make relevant 
conjectures.  We first construct numerically solitary waves, then discuss 
their stability and their appearence in the long-time behavior of 
solutions for general initial data. 

\subsection{Solitary waves}
To study solitary waves of the Benjamin equation, i.e., localized 
solutions of (\ref{SWBOBmod}),
\begin{equation}
	-cQ+\frac{1}{2}Q^{2}- \alpha HQ_{x}-\beta Q_{xx}=0,
	\label{bensol}
\end{equation}
we first construct 
them numerically. Note that $Q$ depends on three parameters, $c$, 
$\alpha$, $\beta$, and only one of them can be scaled out by 
rescaling $Q$ and $x$ appropriately. Therefore we will consider the 
dependence on some of these parameters whilst the others are kept 
fixed. 

The numerical approach is as in \cite{KS15} a Fourier method 
with a Newton-Krylov iteration. The idea is to study (\ref{bensol}) 
in Fourier space 
\begin{equation}
	\mathcal{F}(\hat{Q}):=\mathcal{L}\hat{Q}+\frac{1}{2}\widehat{Q^{2}}=0,
	\label{bensolf}
\end{equation}
where $\mathcal{L}=\beta k^{2}-\alpha|k|-c$. The Fourier transform 
will be approximated in standard way with a discrete Fourier 
transform (DFT); for simplicity both are denoted with the same 
symbol. 
This means we approximate a situation on the real line by one on the 
torus with period $2\pi L$, i.e., we choose $x\in 
L[-\pi,\pi]$ with $L>0$ large enough such that the solutions decay to 
machine precision (here $10^{-16}$) at the boundaries of the 
computational domain. Below we always use $L=20$ for the solitary waves. On the 
interval we use the standard discretisation of the DFT, which is computed via a Fast Fourier transform (FFT). This 
leads for (\ref{bensolf}) to a system of $N$ (the number of discrete 
Fourier modes) nonlinear equations $\mathcal{F}=0$ which are solved 
with a Newton iteration,
$$\hat{Q}^{(n+1)}=\hat{Q}^{(n)}-\left(\mbox{Jac}(\mathcal{F})|_{\hat{Q}^{(n)}}\right)^{-1}\mathcal{F}(\hat{Q}^{(n)}),
\quad n=0,1,\ldots,$$
where $\mbox{Jac}(\mathcal{F})$ is the Jacobian of $\mathcal{F}$. 
To compute 
the action of the Jacobian on a vector, we use the Krylov subspace 
technique GMRES. 

\begin{figure}[!htb]
\includegraphics[width=0.7\hsize]{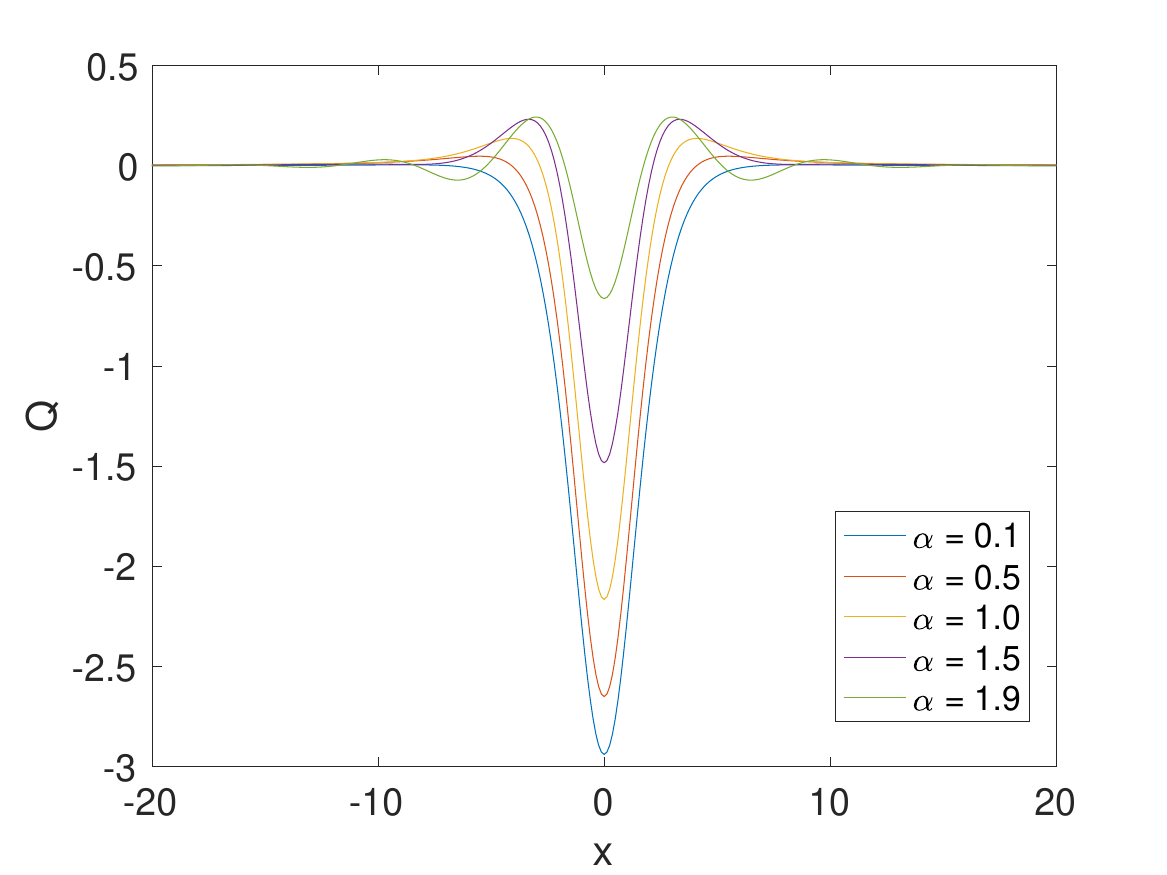}
\caption{Benjamin soliton for $c=-1$, $\beta=1$ and various values of $\alpha$. }
\label{benjaminsol}
\end{figure}

To obtain the solitary waves, we apply a tracing technique. We start 
with the well known KdV 
soliton for $\beta=1$, $c=-1$ and $\alpha=0$, 
$Q=3c\mbox{sech}(x/2/\sqrt(-c/\beta))^2$. We do not change the value of 
$\beta$ and $c$ in the following, but slowly increase the value of $\alpha$ up 
to close to the limiting value 2. The equation for the soliton 
with $c=-1$ is then solved for a slightly larger value of $\alpha$ 
with the solitary wave for the previous value of $\alpha$ as an 
initial iterate $\hat{Q}^{(0)}$.  Note that we always use $N=2^{10}$ Fourier 
modes here, and that the coefficients of the DFT decrease for all 
shown examples to $10^{-10}$ and below which gives an indication of 
the numerical error. We show the resulting solitons for several values of 
$\alpha$ in Fig.~\ref{benjaminsol}. It can be seen that for 
larger values of $\alpha$, the soliton develops some sort of 
elevation which turns into additional oscillations. 

\begin{figure}[!htb]
\includegraphics[width=0.49\hsize]{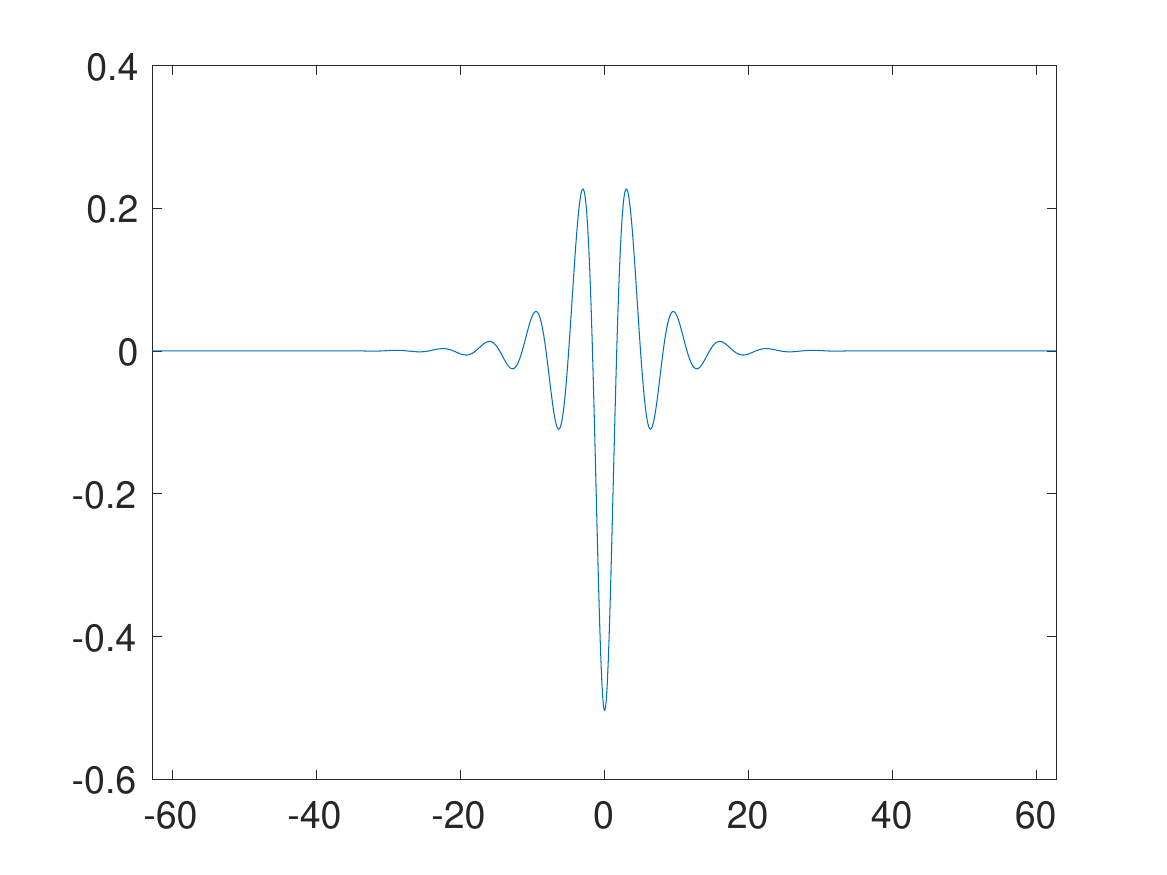}
\includegraphics[width=0.49\hsize]{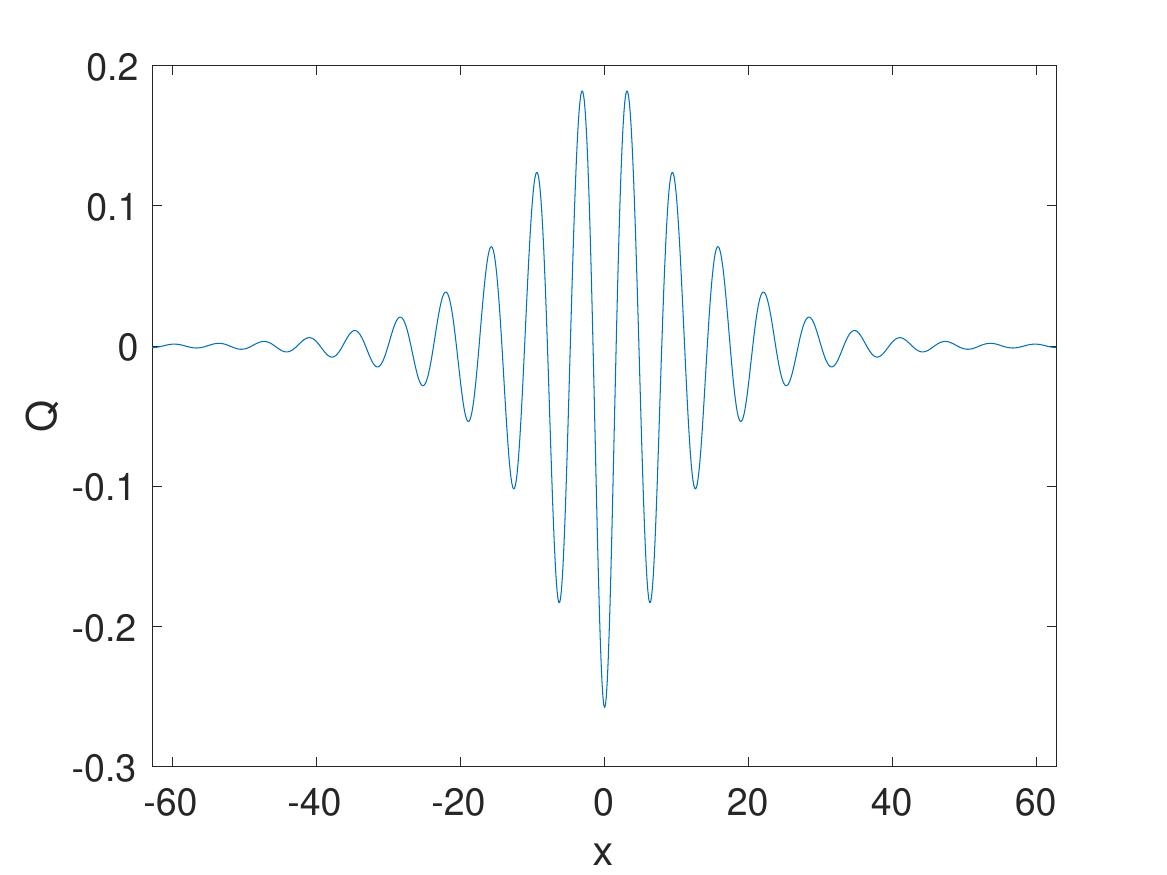}
\caption{Benjamin soliton for $c=-1$, $\beta=1$ and $\alpha=1.95$ on 
the left and $\alpha=1.99$ on the right. }
\label{benjaminsola195}
\end{figure}

For even larger values of $\alpha$, more and more oscillations appear in the limit $\alpha\to2$, see 
Fig.~\ref{benjaminsola195}. Note also that the amplitude of the 
solitons decreases with $\alpha$. 

\begin{figure}[!htb]
\includegraphics[width=0.49\hsize]{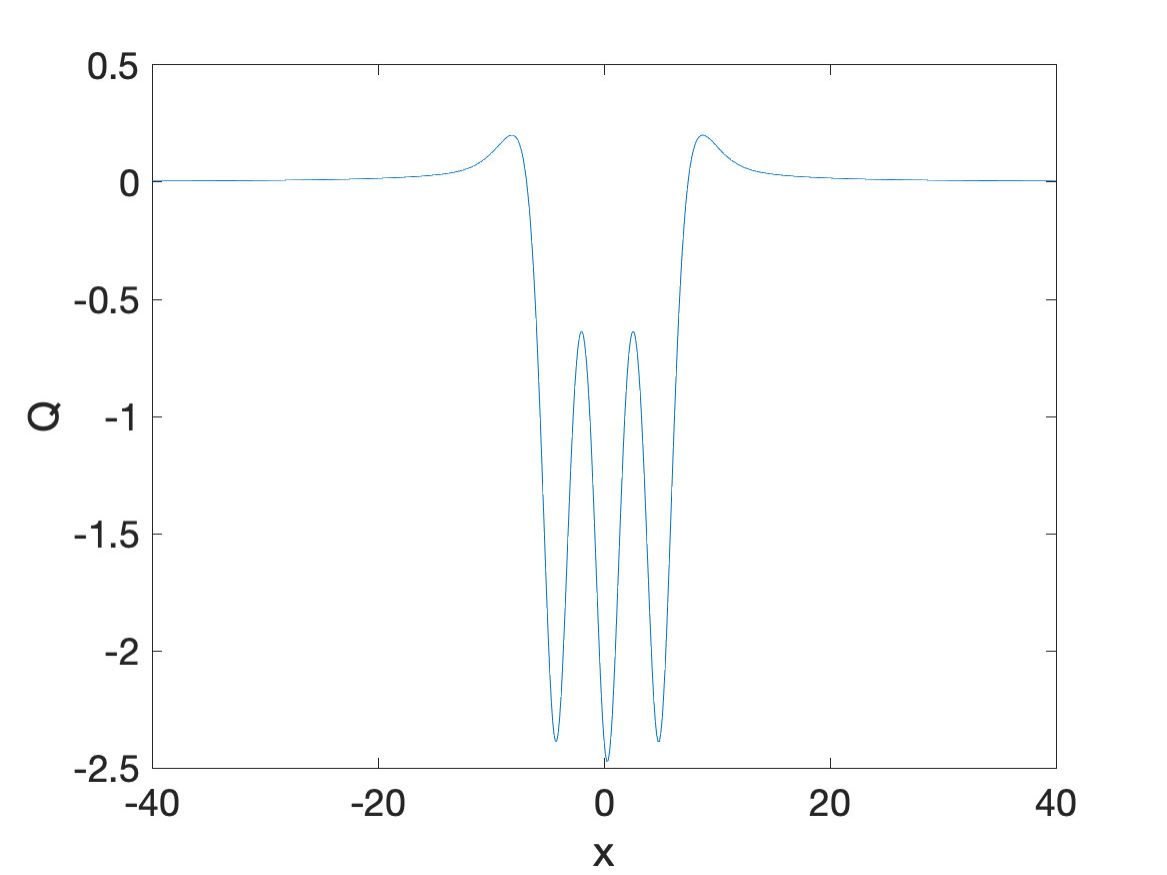}
\includegraphics[width=0.49\hsize]{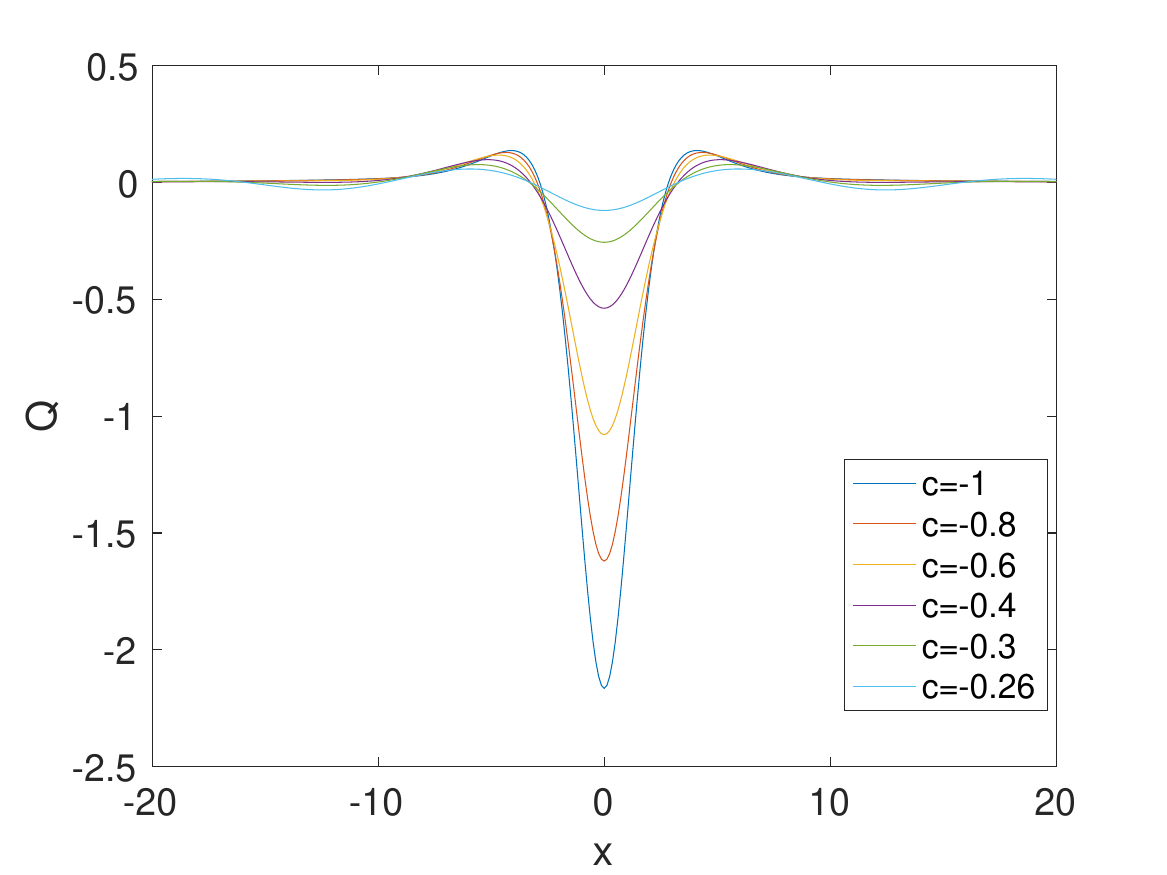}
\caption{On the left the solution to (\ref{bensol}) for $c=-1$, 
$\beta=1$ and $\alpha=1$ for the initial iterate 
$Q=3c\mbox{sech}^2(x/2/\sqrt{-c/\beta})$ (KdV soliton), on the right 
solitary waves for $\alpha=1$, $\beta=1$ and various values of $c$. }
\label{benjaminsola1b1}
\end{figure}

Note that the tracing approach is important to find the solutions 
shown in Fig.~\ref{benjaminsol}. If we do for instance the same 
iteration as above for $\alpha=\beta=1$ with the KdV soliton as 
the initial iterate, we get the solution shown in 
Fig.~\ref{benjaminsola1b1} on the left. The residual of the Newton 
iteration, i.e., $||\mathcal{F}(\hat{Q})||_{\infty}$ 
is as for the solutions in Fig.~\ref{benjaminsol} smaller 
than $10^{-10}$. This could indicate that the critical points  are not 
unique. But if we compute the energy, it is roughly 3.5 times the 
energy of the solitary wave in Fig.~\ref{benjaminsol}. Thus the 
ground states appear to be of the form shown in 
Fig.~\ref{benjaminsol}.  

If we fix $\alpha=\beta=1$ and vary $c$ starting from $c=-1$ to 
larger values of $c$ --- again taking the solution at the slightly 
smaller value of $c$ as the initial iterate for a larger value of 
$c$ --- we get the sequence of solutions shown in 
Fig.~\ref{benjaminsola1b1} on the right. Again the solutions become 
smaller and more oscillatory when $c$ approaches the critical 
velocity $-1/4$. 

\begin{figure}[!htb]
\includegraphics[width=0.49\hsize]{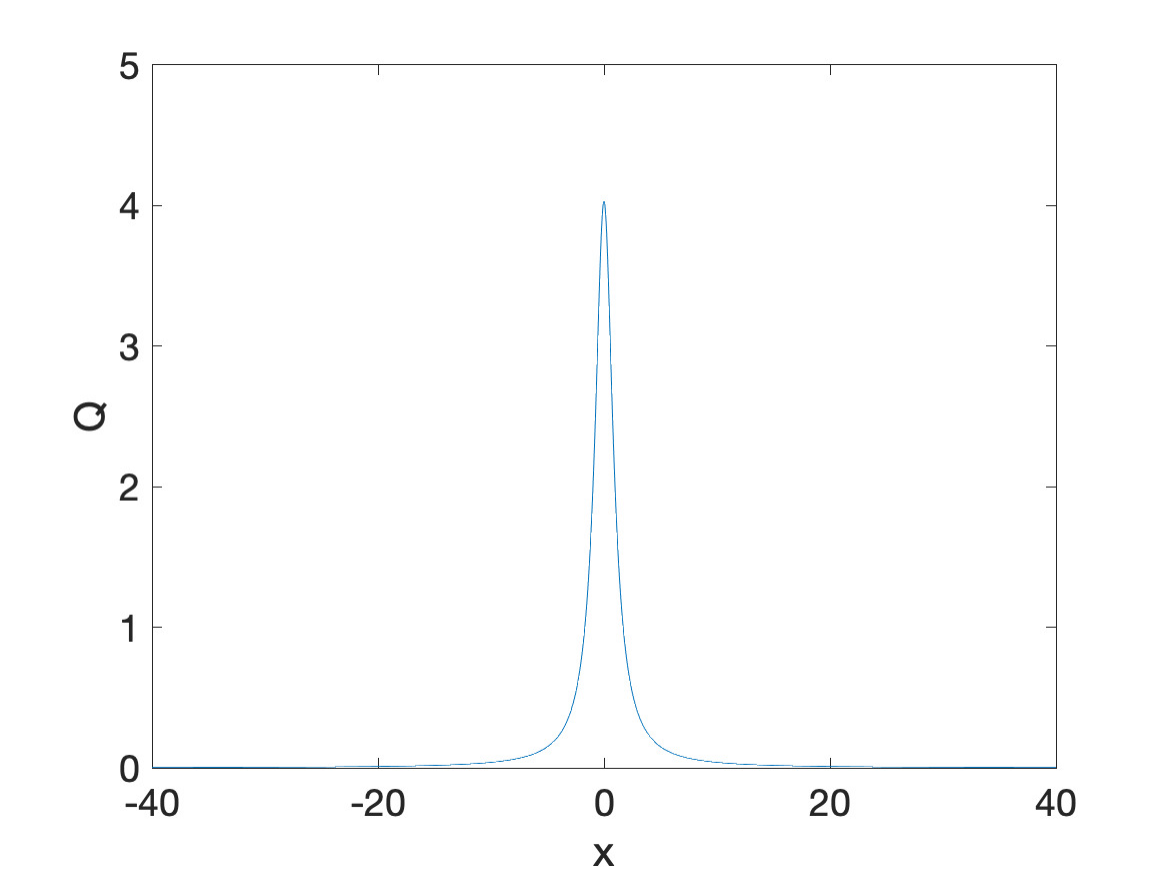}
\includegraphics[width=0.49\hsize]{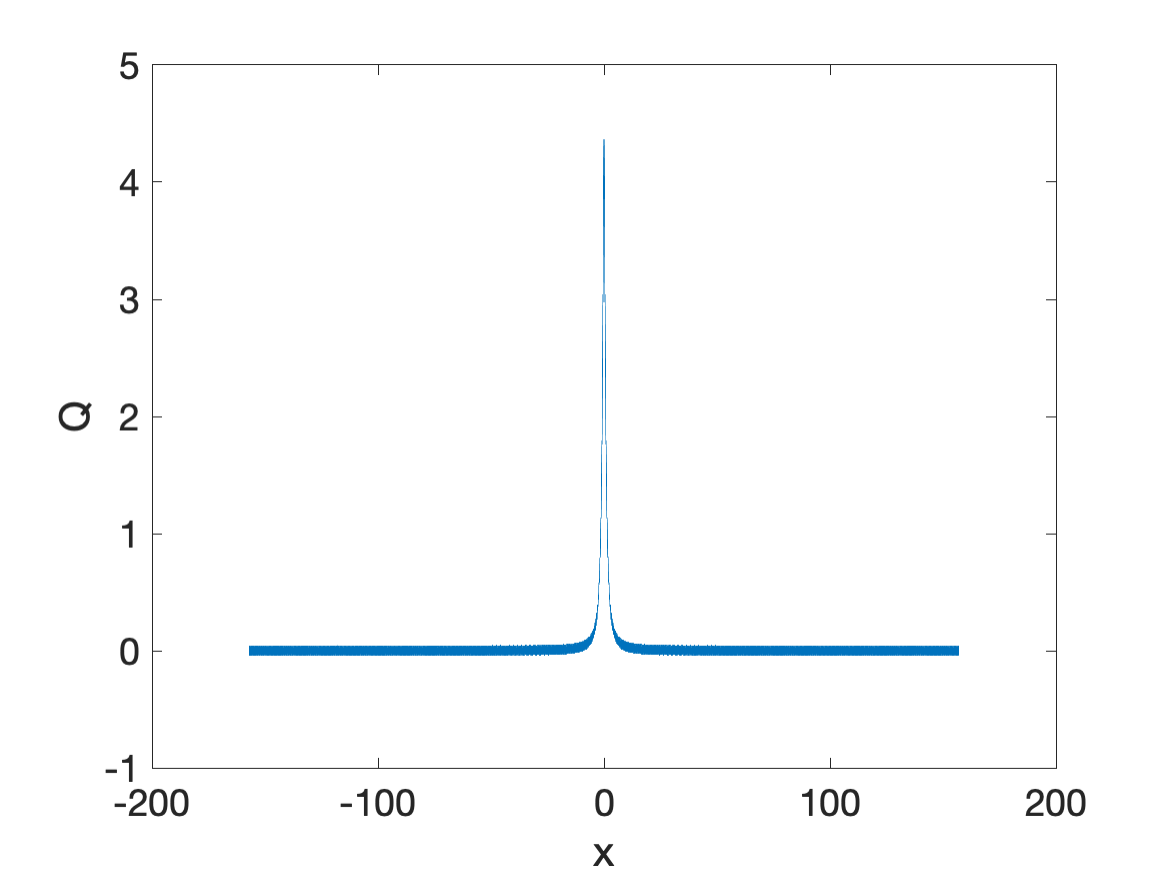}
\caption{Solution to (\ref{bensol}) for $c=1$, $\alpha=1$ and
$\beta=10^{-2}$ on the left and for $c=1$, $\alpha=1$ and
$\beta=6*10^{-2}$ on the right. }
\label{benjaminc1sola1b1em1}
\end{figure}

The situation is slightly different, if we look for solitary waves in 
the vicinity of the BO soliton, $Q = 4c/(1+c^2x^2)$ for $\alpha=1$, 
$\beta=0$. Using the BO soliton for $c=\alpha=1$ and $\beta=10^{-2}$ as 
the initial iterate with $N=2^{12}$ and $L=50$, we find the solution on the left of 
Fig.~\ref{benjaminc1sola1b1em1}. It appears that there are solitary 
waves of a similar form as for BO for small $\beta$. The maximum is 
only slightly larger than the one of the BO soliton. 
However, this changes for values of $\beta\gtrsim5*10^{-2}$. As 
shown on the right of Fig.~\ref{benjaminc1sola1b1em1}, there appears 
to be only a solution to (\ref{bensol}) with an oscillatory 
singularity at infinity. This is not due to the chosen numerical 
approach since the behavior does not depend strongly on the chosen 
period. Even if we apply the method of \cite{hilbert} where equation 
(\ref{bensol}) is considered on the compactified real line without 
any periodicity assumption, these 
oscillations are observed (note that in this case a first order zero 
of $Q$ at infinity is enforced). Thus it appears that solitary waves 
only exist for very small $\beta$ in the vicinity of the BO soliton. 

\begin{remark}
	The iterative approach applied in this section implies that an 
	initial iterate in the vicinity of the wanted solution has to 
	be chosen since the convergence of Newton iterations is local, 
	see for instance Fig.~\ref{benjaminsola1b1} on the left. 
	Consequently we cannot decide with this method whether there are 
	stationary solutions to the Benjamin equation since we do not 
	have an appropriate initial iterate and since the branches of 
	solutions we could study above do not have a (non-trivial) limit 
	for $c\to0$.  
\end{remark}

\begin{figure}[!htb]
\includegraphics[width=0.7\hsize]{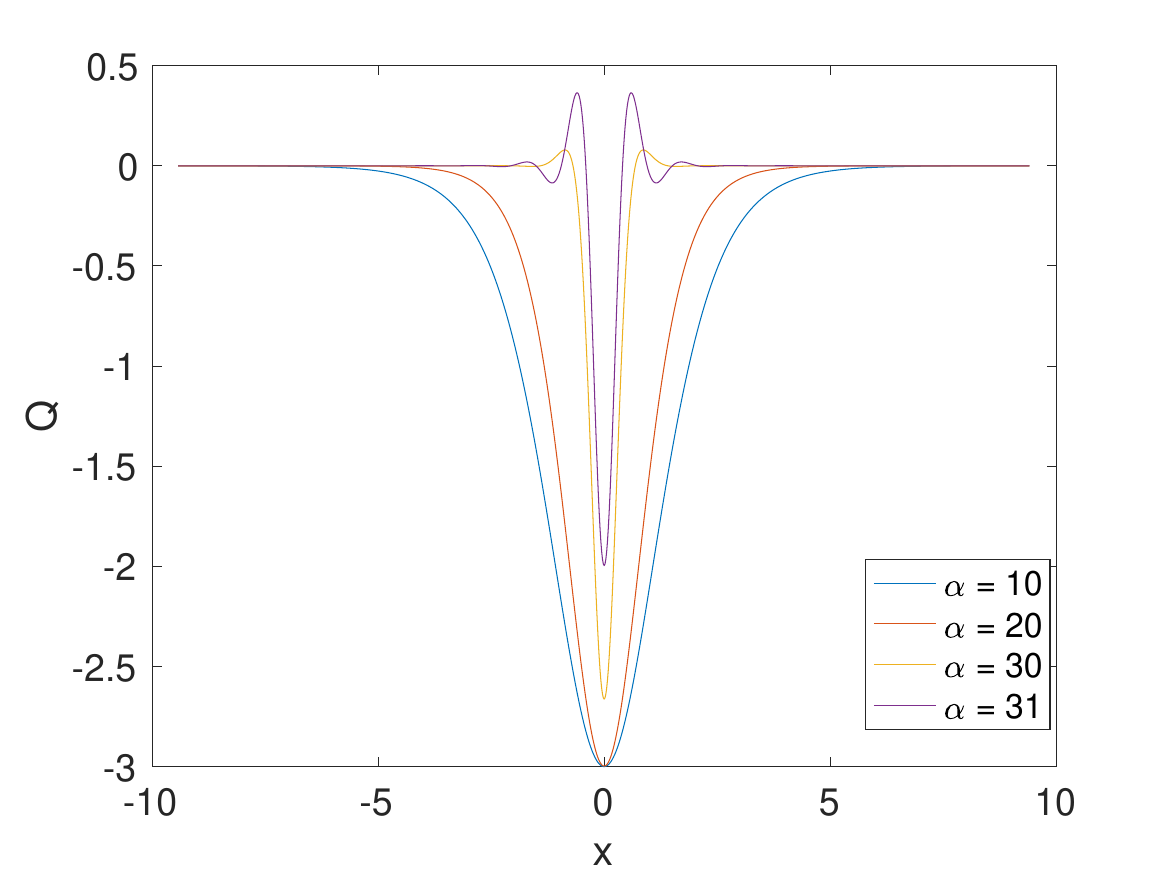}
\caption{Solitary waves for the KdV-ILW equation (\ref{ILWTS}) for 
$c=-1$, $\beta=1$, $\delta=0.1$ and various values of $\alpha$. }
\label{BenIWLd01solcm1alpha}
\end{figure}

The same numerical approach will be used for solitary waves of the KdV-ILW equation 
(\ref{ILWTS}). In Fig.~\ref{BenIWLd01solcm1alpha} we show solitary 
waves obtained by deforming the KdV soliton as in 
Fig.~\ref{benjaminsol}. We put $c=-1$, $\beta=1$, $N=2^{10}$ for 
$x\in3[-\pi,\pi]$, and $\delta=0.1$. The solitary waves are much more localised than 
the ones with a BO term in the Benjamin equation in 
Fig.~\ref{benjaminsol}.

\begin{figure}[!htb]
\includegraphics[width=0.49\hsize]{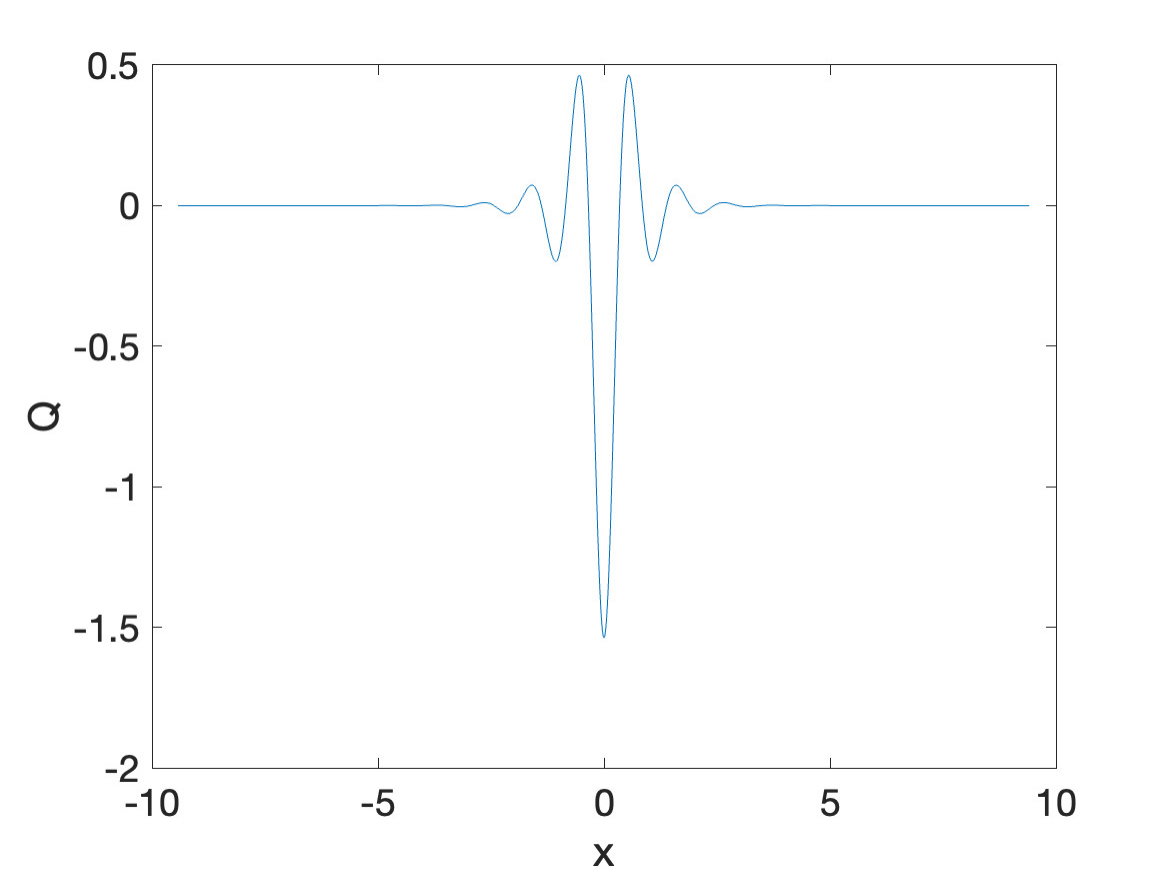}
\includegraphics[width=0.49\hsize]{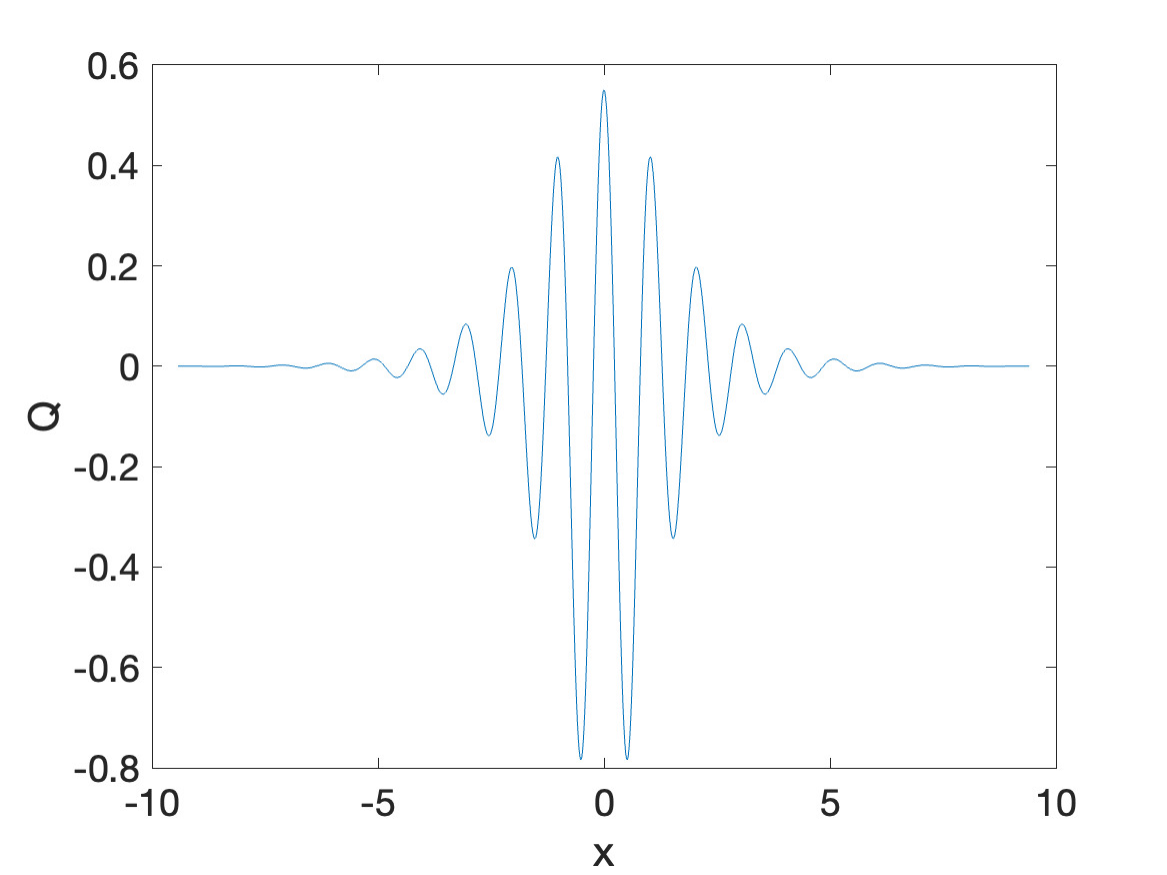}
\caption{Solitary waves for the KdV-ILW equation (\ref{ILWTS}) for 
$c=-1$, $\beta=1$, $\delta=0.1$ and $\alpha=31.3$ on 
the left and $\alpha=31.5$ on the right. }
\label{BenIWLd01solcm1alpha2}
\end{figure}

Again the waves become smaller in amplitude and more oscillatory for 
larger values of $\alpha$, which is even more visible in 
Fig.~\ref{BenIWLd01solcm1alpha2}. The limiting value of $\alpha$ 
beyond which there are no more localised solitary waves in this 
sequence is unknown, but must be close to the shown values.

\begin{figure}[!htb]
\includegraphics[width=0.49\hsize]{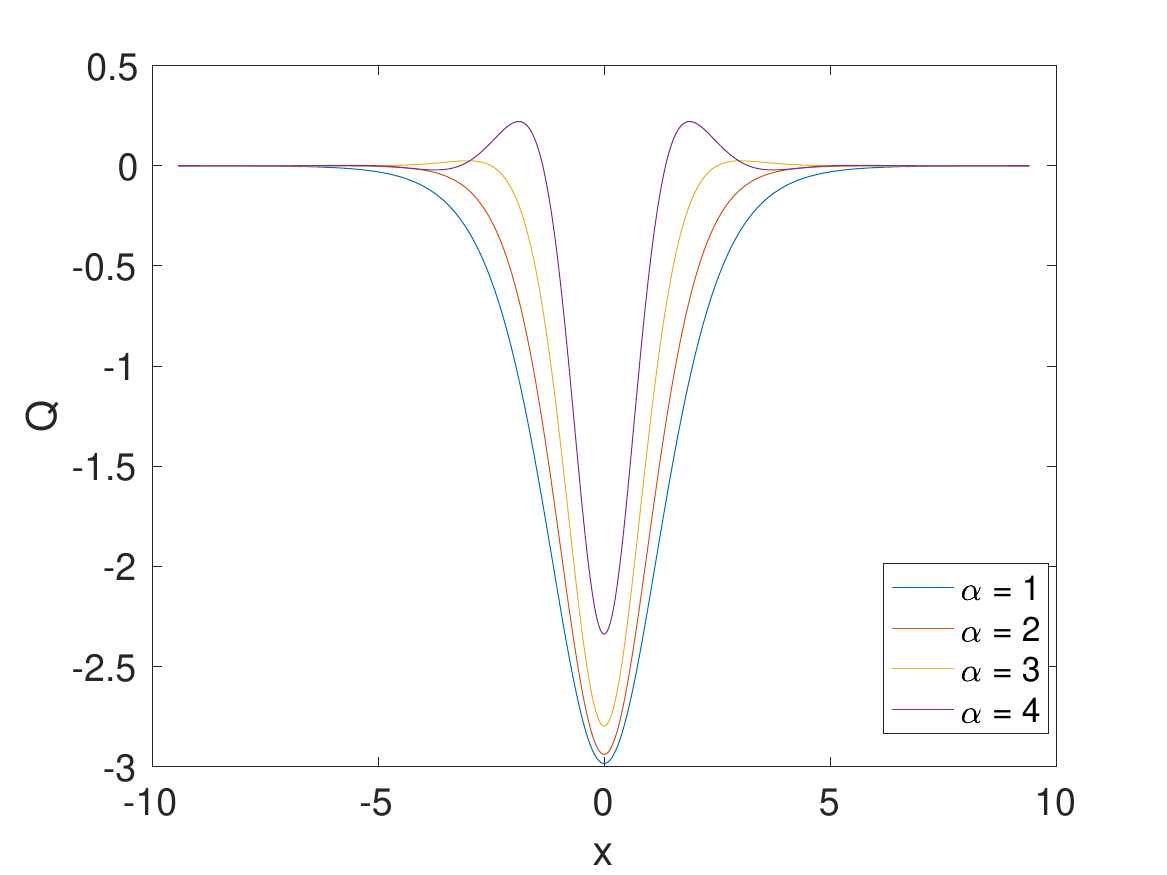}
\includegraphics[width=0.49\hsize]{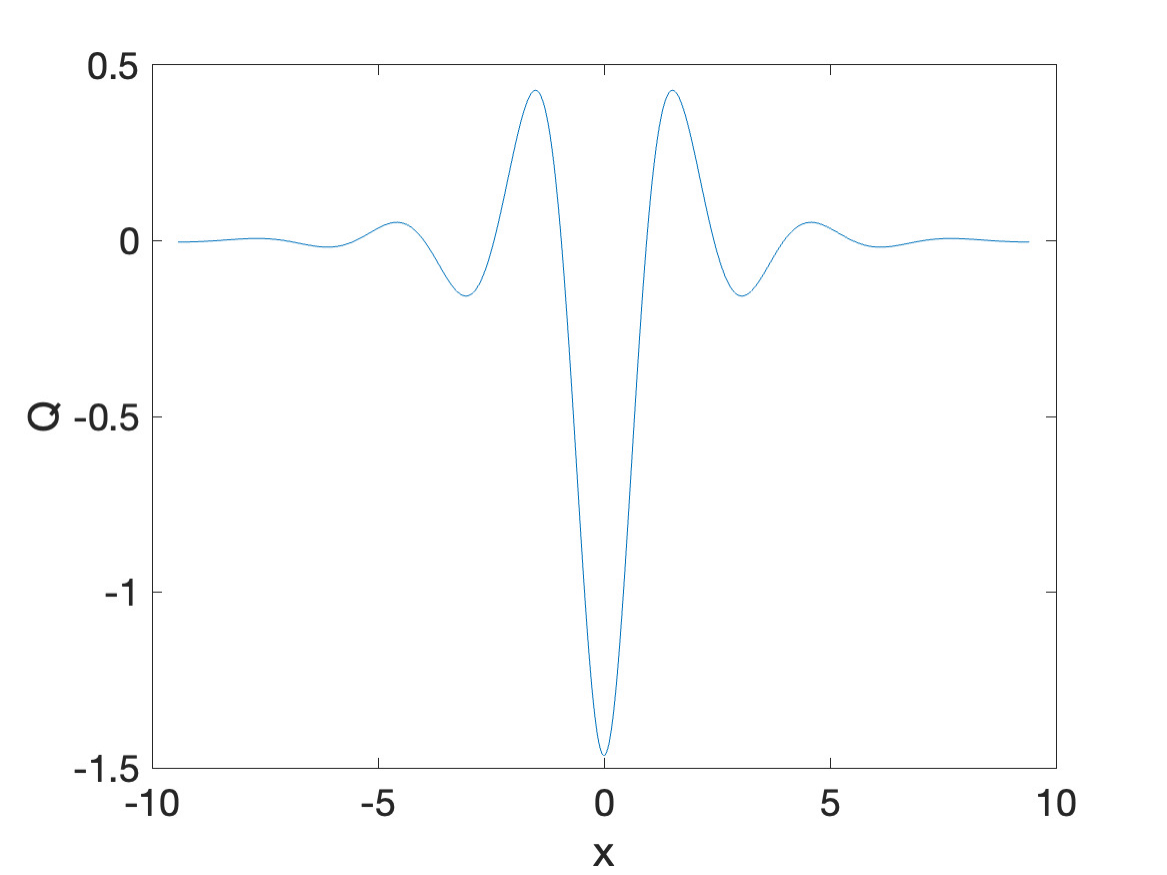}
\caption{Solitary waves for the KdV-ILW equation (\ref{ILWTS}) for 
$c=-1$, $\beta=1$, $\delta=0.9$ and several values of $\alpha$ on 
the left and for $\alpha=4.7$ on the right. }
\label{BenIWLd09solcm1alpha}
\end{figure}

For larger values of $\delta$, here $\delta=0.9$, the behavior is 
similar, but closer to the KdV soliton as can be seen in 
Fig.~\ref{BenIWLd09solcm1alpha}. The solitary wave becomes once more 
oscillatory for $\alpha=4.7$ as shown on the right of the same 
figure. The precise limiting value of $\alpha$ for solitary waves in 
this sequence is again unknown but close to 5.  

\begin{figure}[!htb]
\includegraphics[width=0.49\hsize]{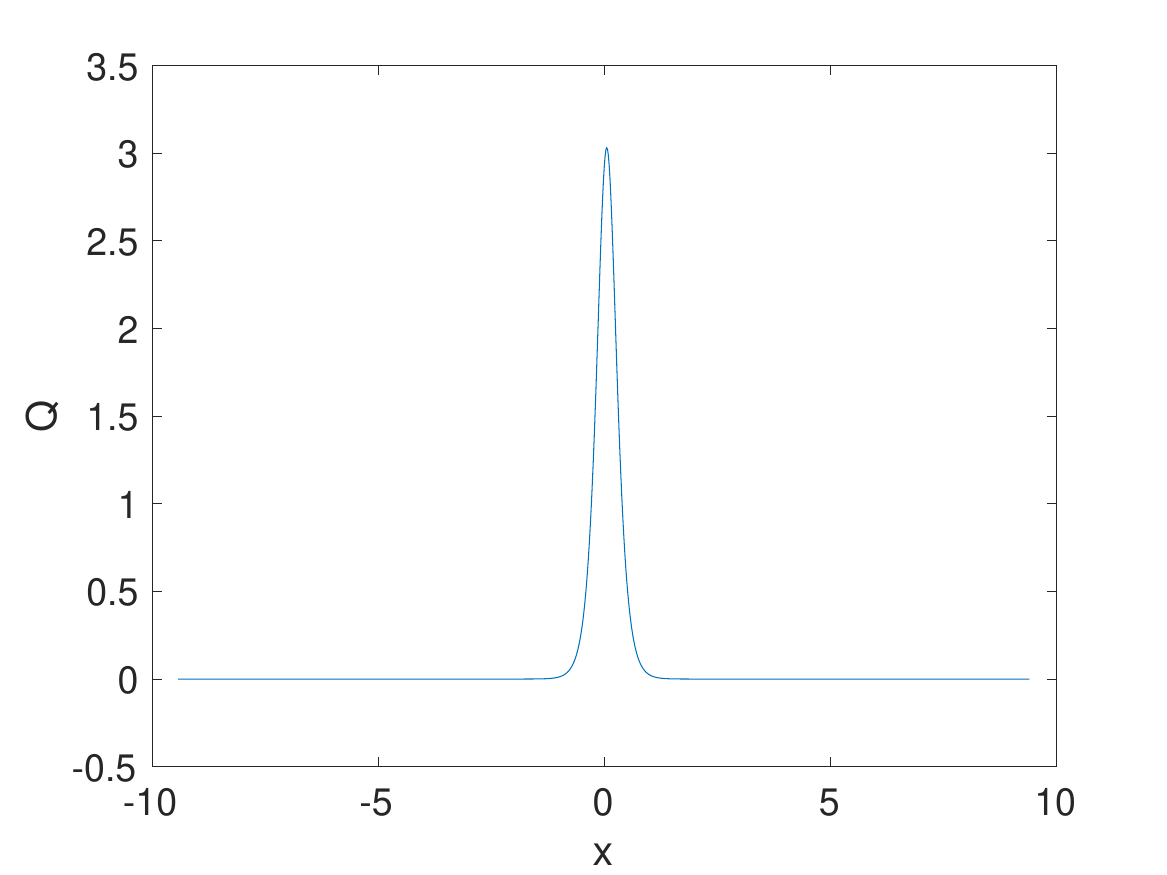}
\includegraphics[width=0.49\hsize]{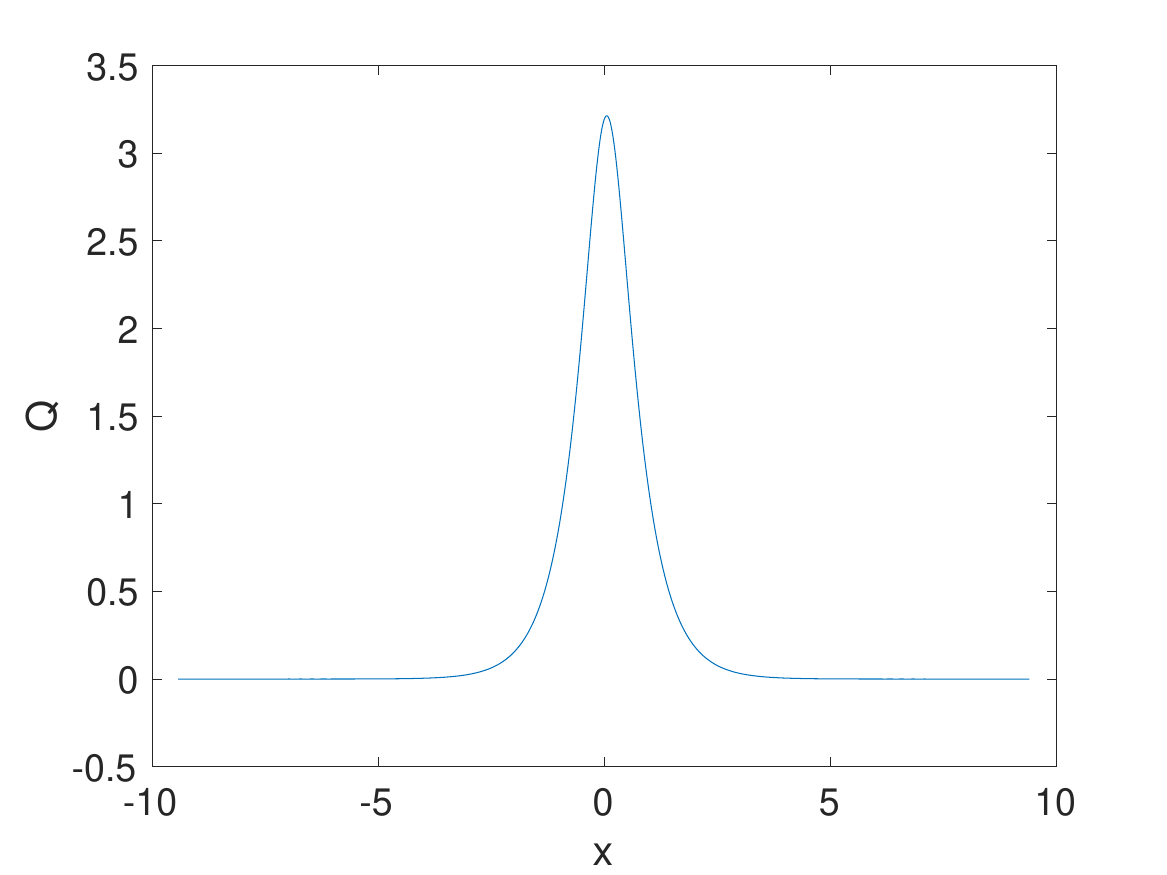}
\caption{Solitary waves to KdV-ILW equation to (\ref{ILWTS}) for 
$c=1$, $\alpha=1$, $\delta=0.1$ and
$\beta=10^{-2}$ on the left and for $c=1$, $\alpha=1$, $\delta=0.9$ and
$\beta=4*10^{-2}$ on the right. }
\label{BenIWLd09solc1a1dela}
\end{figure}

The ILW equation has explicit solitons where only parameters have to 
be determined for given velocity $c$ and $\delta$ as the solution of 
a transcendental equation, see for instance \cite{KS}. As in Fig.~\ref{benjaminc1sola1b1em1} one 
can thus look for solitary waves in the vicinity of the ILW soliton. 
We put $c=1$, $\alpha=1$ and consider small values of $\beta$. We get 
solitary waves for small values of $\beta$ in this case that are 
close to the ILW solitons, see Fig.~\ref{BenIWLd09solc1a1dela}. 
For larger values of $\beta$ than those 
shown in the figure, rapid oscillations appear quickly as in 
Fig.~\ref{benjaminc1sola1b1em1} on the right. 

\subsection{Time evolution}

For the time evolution of initial data for the Benjamin equation, we 
use the same discretisation in $x$ and FFT techniques as in the 
construction of the soliton. This leads to a system of the form 
$$\hat{Q}_{t}=\mathbb{L}\hat{Q}+\mathcal{N}(\hat{Q}),$$
where $\mathbb{L}=ik(\alpha|k|-\beta k^{2})$ and where 
$\mathcal{N}(\hat{Q})=ik\widehat{Q^{2}}/2$. Since  $\mathbb{L}$ is 
cubic in $k$, the resulting system is \emph{stiff} which means that 
explicit time stepping schemes are inefficient. Therefore we apply
an exponential time differencing (ETD) scheme, see \cite{HO} for 
a review. Since in \cite{KR} the studied fourth order ETD schemes 
performed similarly, we apply here the scheme by Cox and Matthews 
\cite{CM}. The numerical accuracy of the solution is controlled as 
discussed in \cite{KR} via the decrease of the Fourier 
coefficients\footnote{It is known that the DFT of analytic periodic functions decreases as 
the Fourier coefficients exponentially, and for simplicity we speak 
in the following of the Fourier coefficients when we refer to the 
DFT.} and the conservation of the numerically conserved energy. The 
latter will depend on time due to unavoidable numerical errors, and 
its conservation will thus allow to estimate the numerical accuracy 
of the time integration. 

The code is tested at the example of a solitary wave, for instance 
the one for $\alpha=1.95$, $\beta=1$ and $c=-1$ shown in 
Fig.~\ref{benjaminsola195} on the left. We use $N=2^{10}$ Fourier 
modes for $x\in 20[-\pi,\pi]$ with $N_{t}=10^{4}$ time steps for 
$t\in[0,4]$. The energy is conserved relatively to the order 
$10^{-12}$ during the whole computation. The difference between the 
numerical solution and the initial data propagated with $c=-1$ is of 
the order of $10^{-14}$. This tests both the time evolution code as 
well as the code for the numerical construction of the solitary 
waves. It also shows that the energy conservation can be used as an 
indicator of the acuracy of the time integration. 

\subsubsection{Stability of solitary waves}

\begin{figure}[!htb]
\includegraphics[width=0.49\hsize]{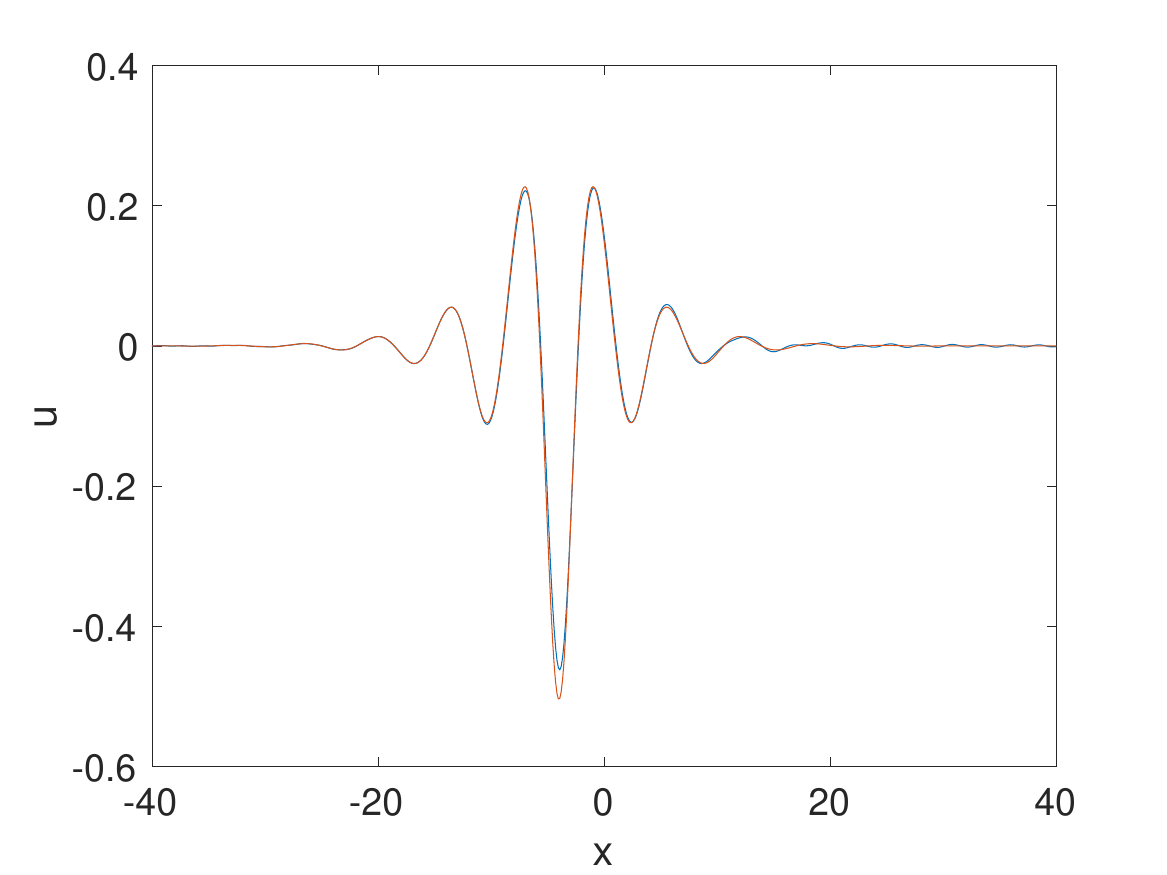}
\includegraphics[width=0.49\hsize]{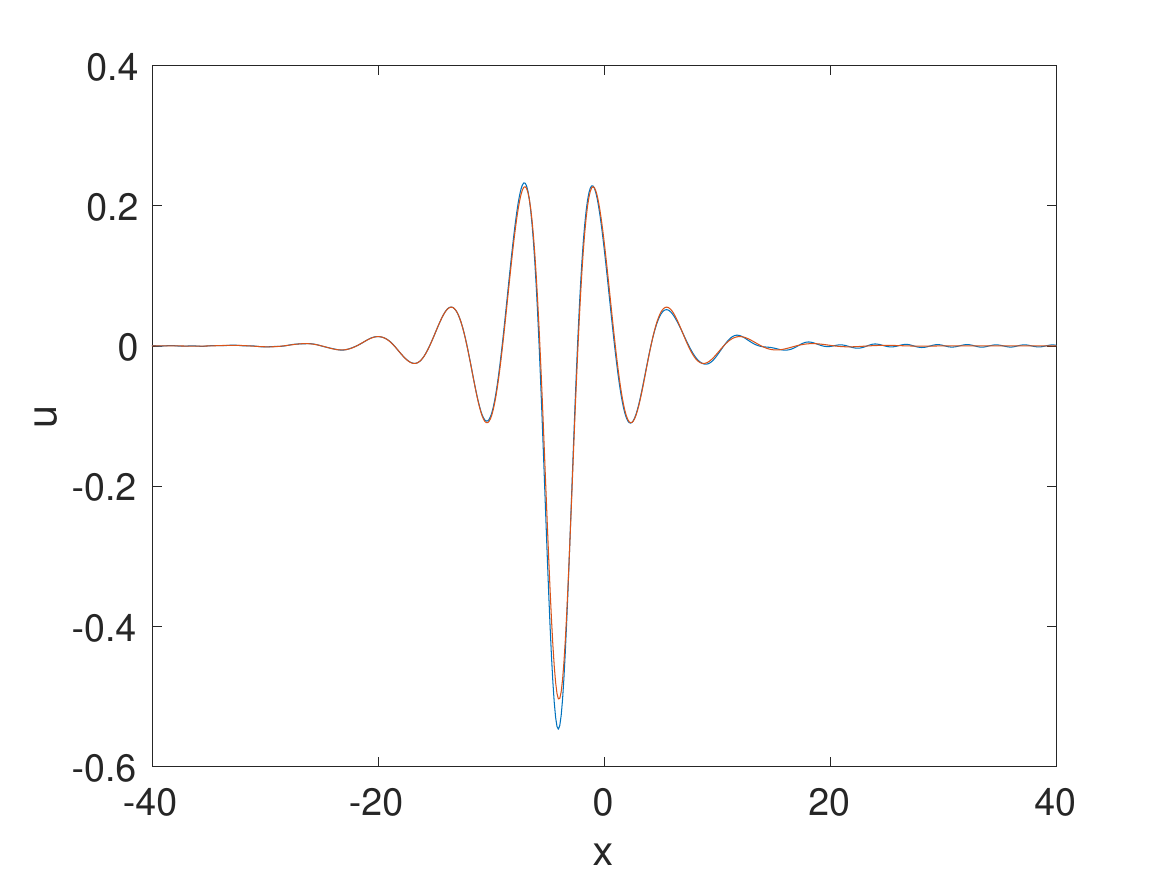}
\caption{Solution to the Benjamin equation for $t=4$ for initial data 
$u_{0}(x)=Q(x)\pm0.05\exp(-x^{2})$ where $Q$ is the solitary wave for 
$\alpha=1.95$, $\beta=1$ and $c=-1$, in blue and the unperturbed 
solitary wave in red; on the left for the $+$ sign of the 
perturbation, on the right for the $-$ sign. }
\label{benb1a195ssolper}
\end{figure}

The Benjamin solitary wave is known to be stable. To illustrate this, we 
consider perturbations of the soliton with $c=-1$, $\alpha=1.95$ and 
$\beta=1$ as above and use the same numerical parameters as there. 
In Fig.~\ref{benjaminsola195}, we show the solution to 
the Benjamin equation for 
$t=4$ for the initial data $u_{0}(x)=Q(x)\pm0.05\exp(-x^{2})$, where $Q$ is the 
solitary wave for these parameters. This means we consider a solitary 
wave with a small Gaussian perturbation.  The numerically computed 
energy is conserved to the order of $10^{-12}$. The solutions at the 
final time are shown in Fig.~\ref{benb1a195ssolper}. In red we give 
the unperturbed solitary wave. Since the perturbation is small, but 
finite (in order to allow for numerically visible effects of the 
perturbation for finite times), the final state is slighly different from the unperturbed 
solitary wave, but obviously as expected close. The perturbation 
leads to radiation emitted towards infinity. 

To explore solitary waves in the vicinity of BO, we consider the case 
$\alpha=c=1$ and $\beta=0.02$. The solitary wave is constructed with 
$N=2^{14}$ Fourier modes for $x=50[-\pi,\pi]$. Initial data 
corresponding to this solitary wave are  
propagated with $N_{t}=10^{4}$ time steps for $t\in[0,4]$ with a 
relative energy of the order of $10^{-12}$ to a solution with a 
difference to the travelling wave of the same order at the final time 
of the computation. 

\begin{figure}[!htb]
\includegraphics[width=0.49\hsize]{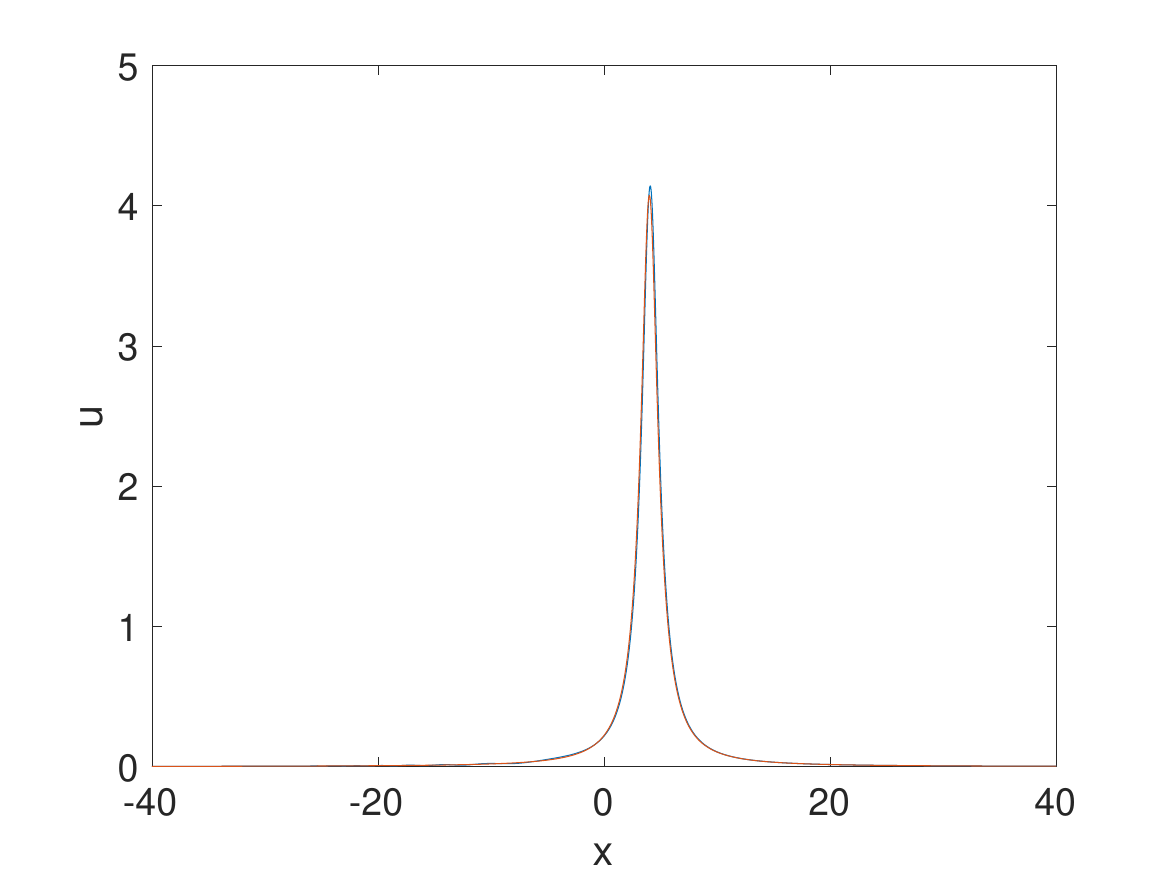}
\includegraphics[width=0.49\hsize]{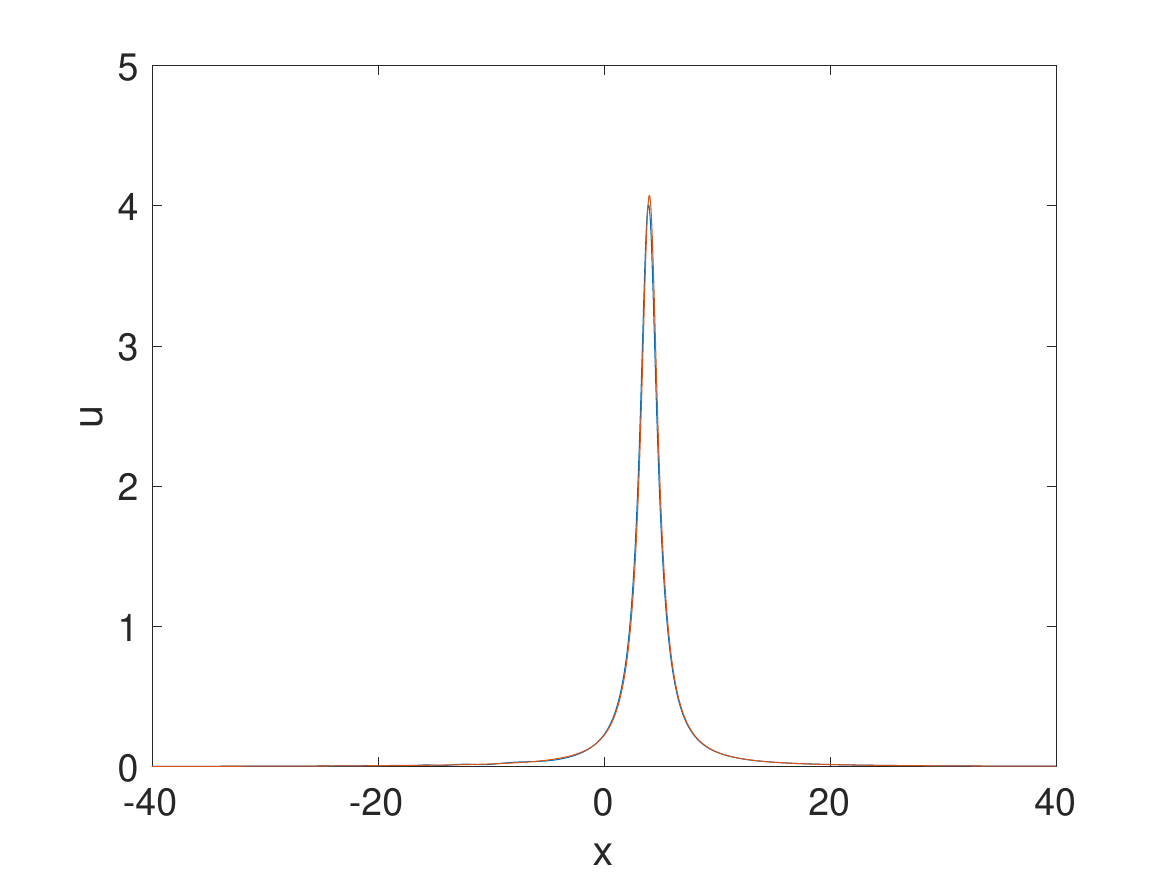}
\caption{Solution to the Benjamin equation for $t=4$ for initial data 
$u_{0}(x)=Q(x)\pm0.05\exp(-x^{2})$ where $Q$ is the solitary wave for 
$\alpha=1$, $\beta=0.02$ and $c=1$, in blue and the unperturbed 
solitary wave in red; on the left for the $+$ sign of the 
perturbation, on the right for the $-$ sign. }
\label{benjaminsolc1a1b2em2per}
\end{figure}

We apply the same numerical parameters for the studies presented 
below. A perturbation of the order of $10\%$ seems to lead to 
instabilities, but these are presumably outside the domain of 
applicability of stability theory. A perturbation of the order of 
$1\%$ leads to similar results as in Fig.~\ref{benb1a195ssolper}. 
Concretely we consider initial data of the form 
$u_{0}(x)=Q(x)\pm0.04\exp(-x^{2})$ where $Q$ is this time the 
solitary wave for $\alpha=c=1$ and $\beta=0.02$. The solutions at time 
$t=4$ are shown in Fig.~\ref{benjaminsolc1a1b2em2per}. They are again 
close to the unperturbed travelling wave which indicates the 
stability of these solitary waves. 

\subsubsection {The soliton resolution conjecture}
\begin{figure}[!htb]
\includegraphics[width=0.7\hsize]{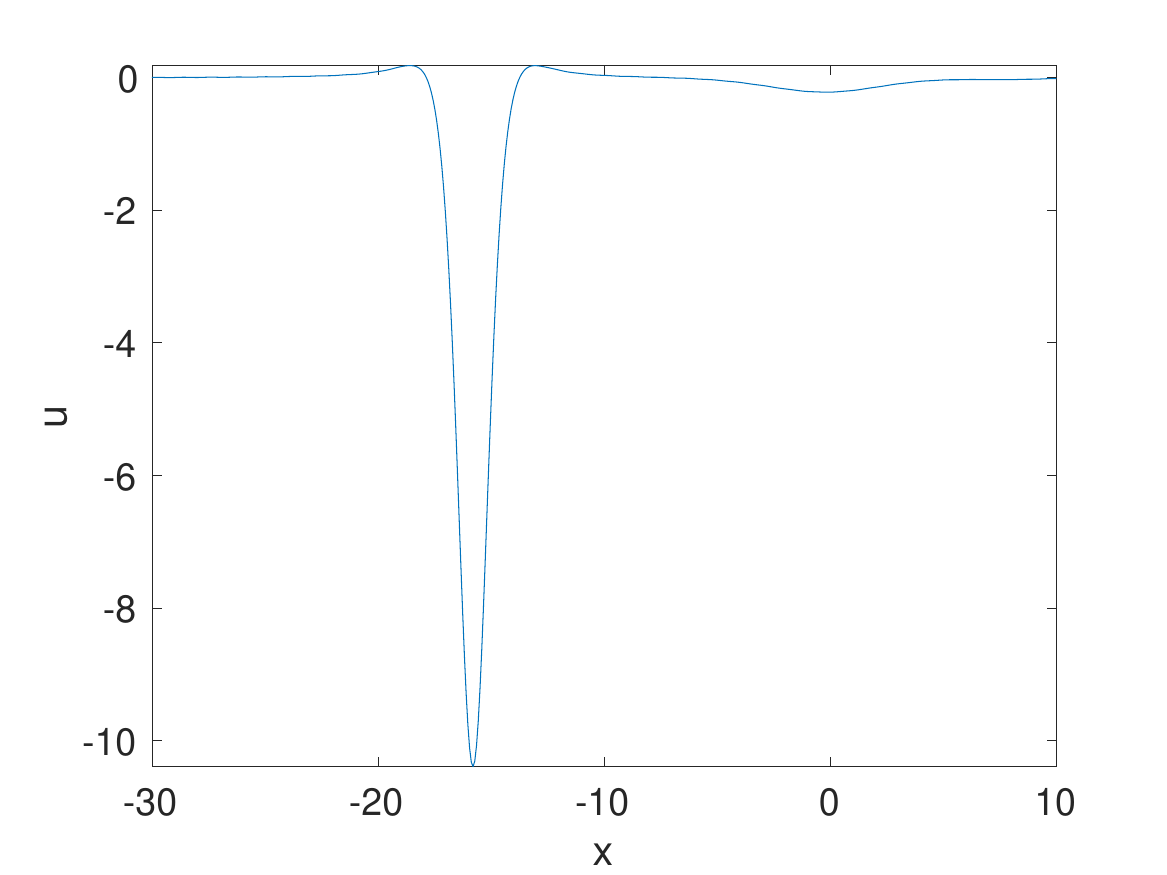}
\caption{A close-up of the solution to the Benjamin equation with $\alpha=\beta=1$ for $t=4$ for initial data 
$u_{0}(x)=-10\exp(-x^{2})$. }
\label{ben10gausszoom}
\end{figure}

Since the solitary waves are stable, it is expected that they appear 
in the long time behavior of solutions of the Benjamin equation for 
initial data of sufficient mass. As an example we consider 
$u_{0}(x)=-10\exp(-x^{2})$ with $N=2^{12}$ Fourier modes and $L=50$ 
(the other numerical parameters are unchanged). The solution for 
$t=4$ can be seen in Fig.~\ref{ben10gausszoom}. In the shown close-up 
at least one soliton 
appears, possibly two. Since the solitary wave does not have a simple 
scaling in $c$, it is difficult to work out the related $c$ and to 
compare with the corresponding solution of (\ref{bensol}). 

\begin{figure}[!htb]
\includegraphics[width=0.49\hsize]{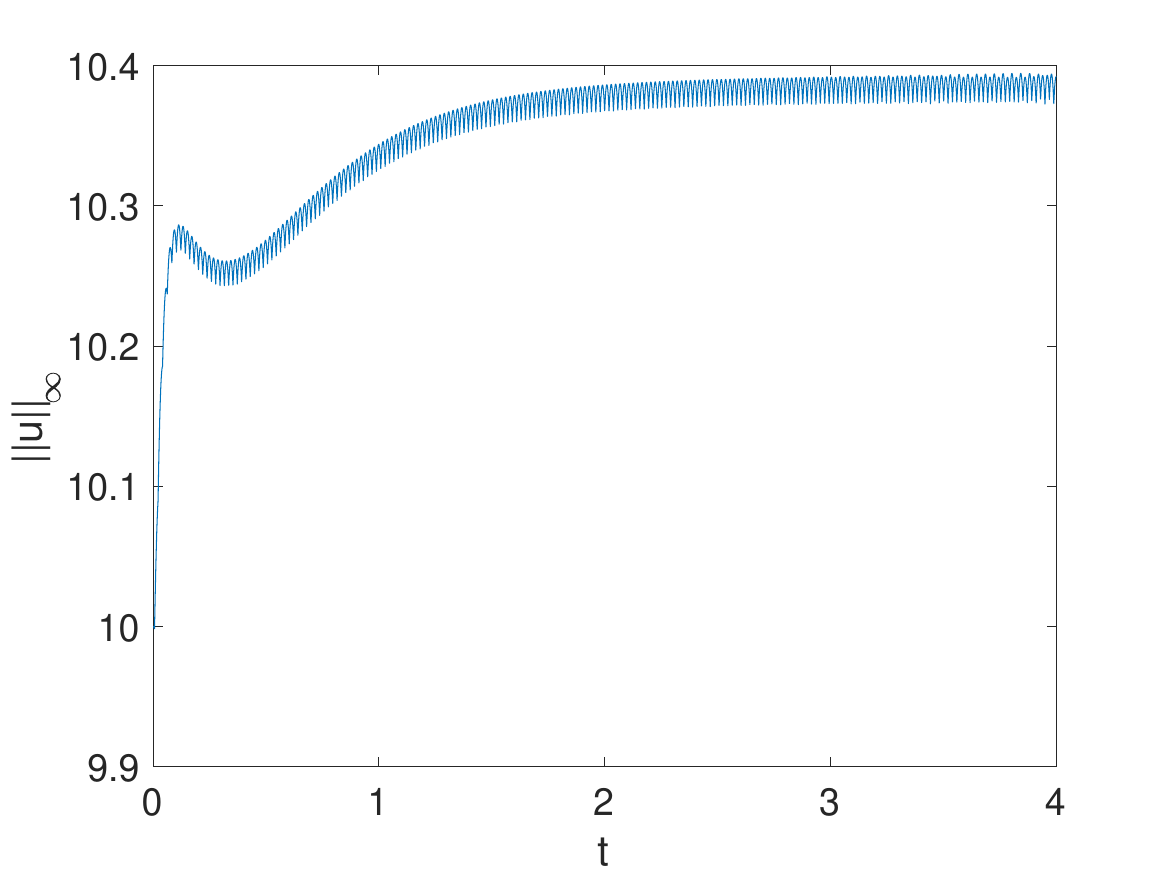}
\includegraphics[width=0.49\hsize]{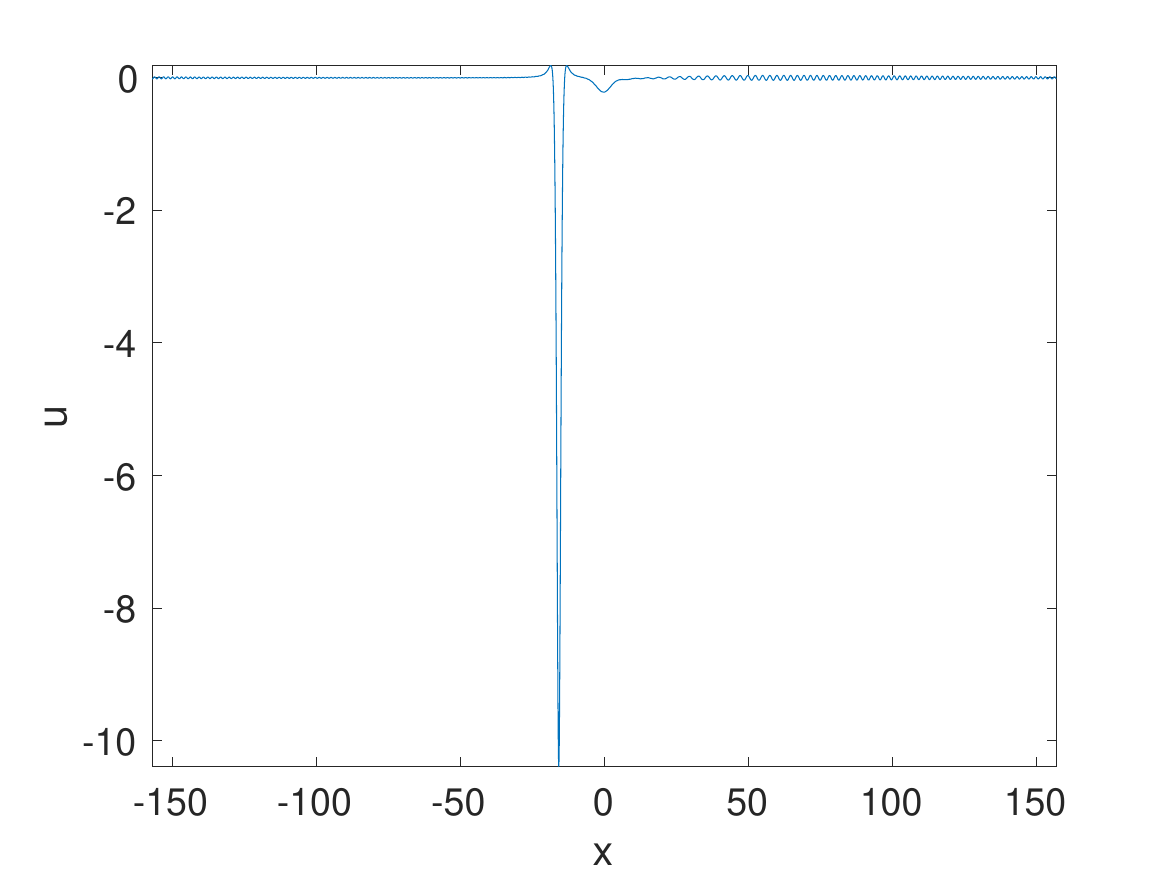}
\caption{$L^{\infty}$ norm of the solution to the Benjamin equation 
with $\alpha=\beta=1$ for initial data 
$u_{0}(x)=-10\exp(-x^{2})$ in dependence of $t$ on the left, and the 
solution of Fig.~\ref{ben10gausszoom} on the full computational 
domain on the right. }
\label{ben10gauss}
\end{figure}

The interpretation of the larger hump as a solitary wave is in 
accordance with  
the $L^{\infty}$-norm of the solution on the left of 
Fig.~\ref{ben10gauss}. It appears to settle on a plateau which is 
generally an indication that a stable solitary wave is the final state 
of the largest hump. The oscillations in the $L^{\infty}$ norm are 
due to radiation emitted to the left which reappears on the other 
side of the computational domain since we approximate a situation on 
$\mathbb{R}$ by a situation on $\mathbb{T}$. The radiation can be 
seen on the right of Fig.~\ref{ben10gauss} where the solution of 
Fig.~\ref{ben10gausszoom} is shown on the full computational 
domain. 

\begin{figure}[!htb]
\includegraphics[width=0.49\hsize]{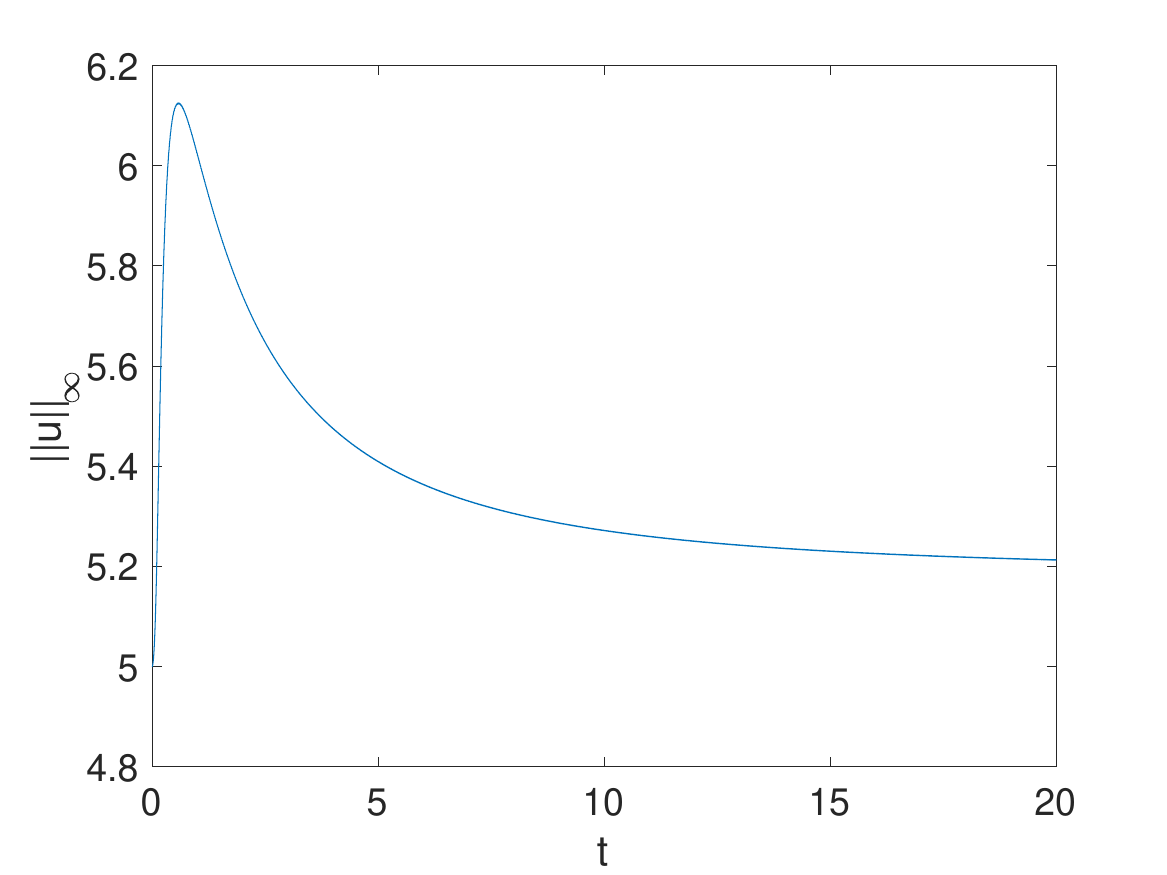}
\includegraphics[width=0.49\hsize]{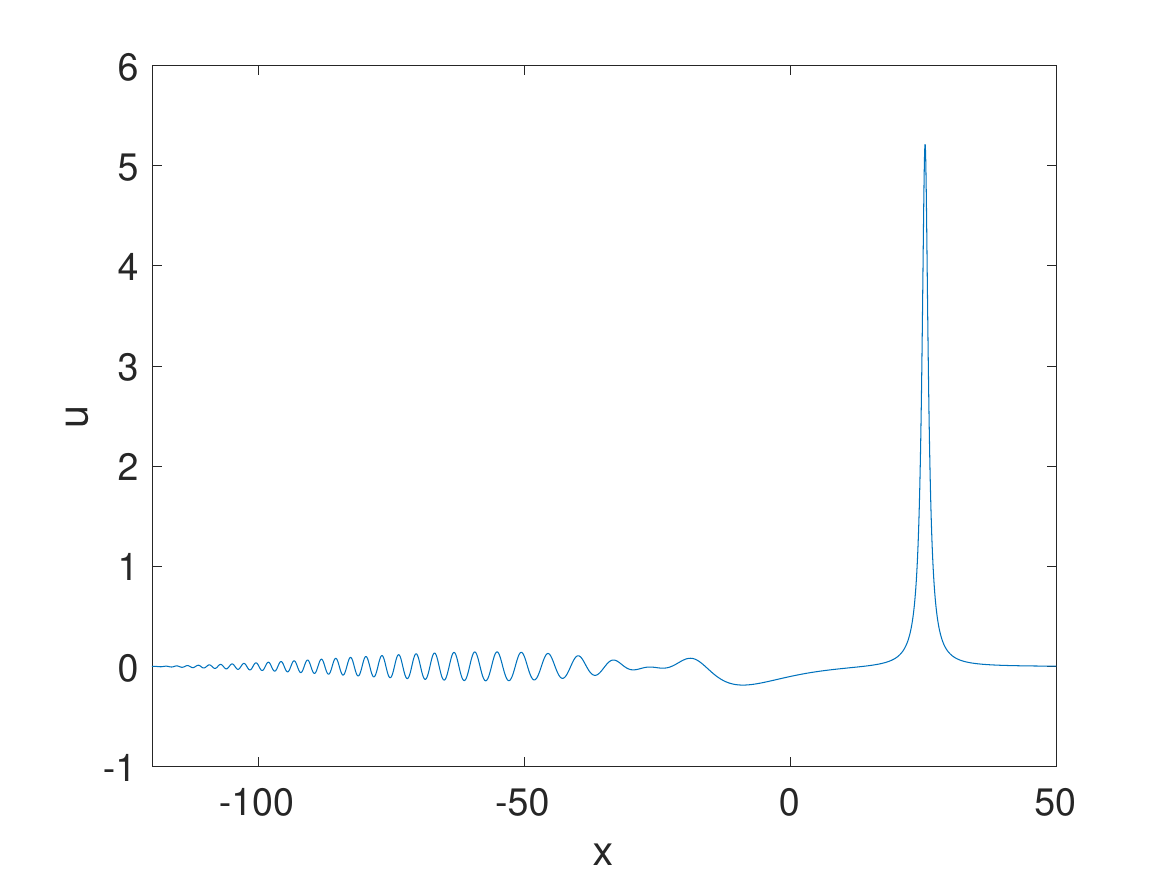}
\caption{$L^{\infty}$ norm of the solution to the Benjamin equation 
with $\alpha=1$, $\beta=0.02$ for initial data 
$u_{0}(x)=5\exp(-x^{2})$ in dependence of $t$ on the left, and the 
solution for $t=20$ on the right. }
\label{ben5gauss}
\end{figure}

The situation is slightly different in the vicinity of BO, for 
instance for $\alpha=1$ and $\beta=0.02$, where we found a solitary 
wave for $c=1$, see Fig.~\ref{benjaminsola1b1} on the left. We use 
$N=2^{14}$ Fourier modes for $x\in50[-\pi,\pi]$ and $N_{t}=10^{4}$ 
time steps for $t\in[0,20]$. The solution for the initial data 
$u_{0}=5\exp(-x^{2})$  can be seen on the right of 
Fig.~\ref{ben5gauss}. There appears to be a solitary wave plus radiation 
traveling to the left. This interpretation is backed by the 
$L^{\infty}$ norm shown on the left of the same figure seemingly 
reaching a plateau for large times. 

\begin{figure}[!htb]
\includegraphics[width=0.49\hsize]{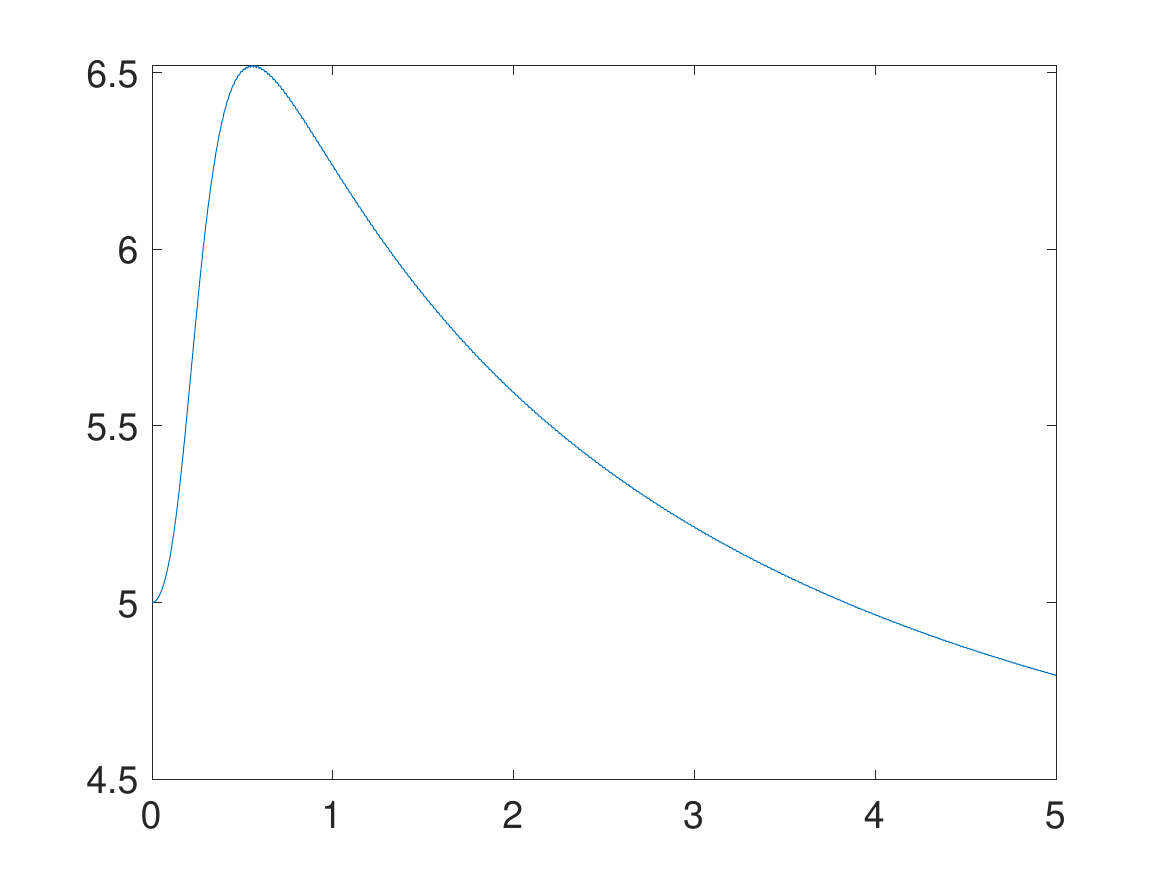}
\includegraphics[width=0.49\hsize]{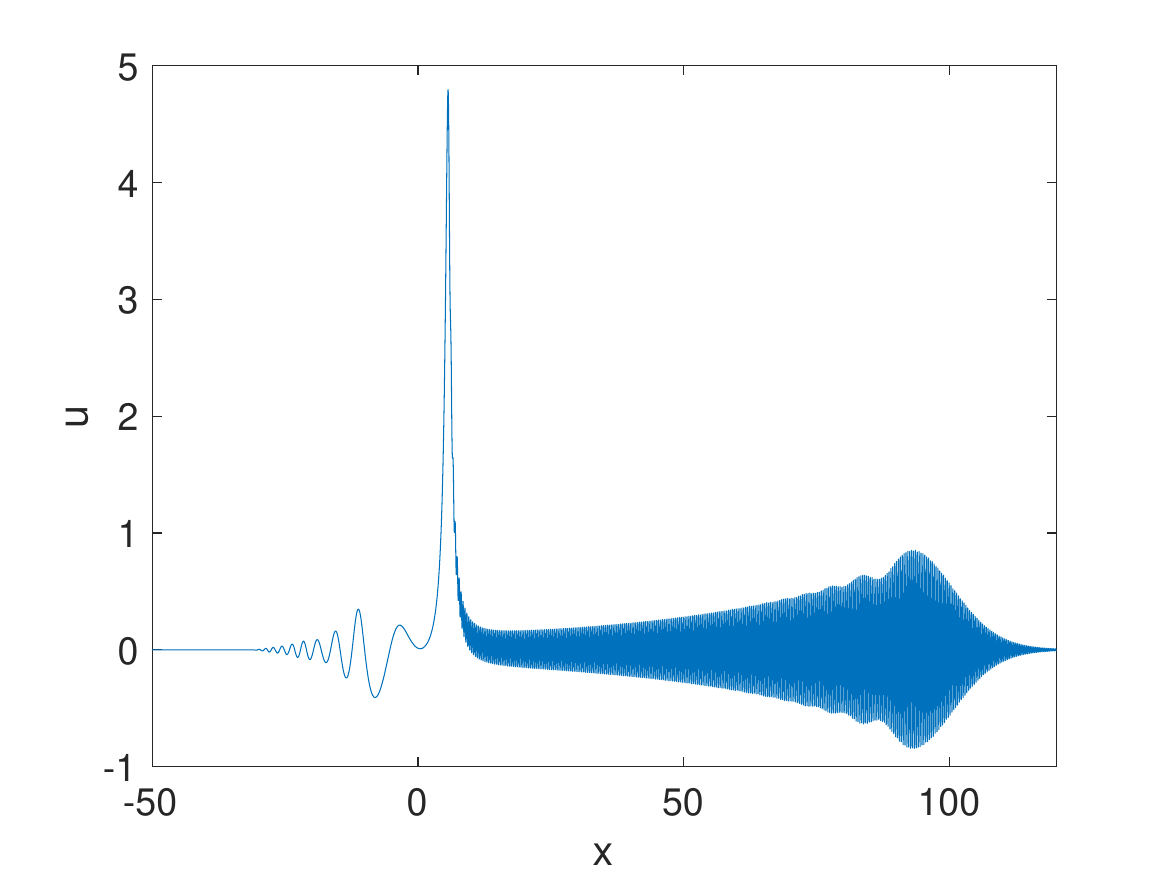}
\caption{$L^{\infty}$ norm of the solution to the Benjamin equation 
with $\alpha=1$, $\beta=0.06$ for initial data 
$u_{0}(x)=5\exp(-x^{2})$ in dependence of $t$ on the left, and the 
solution for $t=5$ on the right. }
\label{ben5gauss6em2}
\end{figure}

If we apply the same numerical parameters for times smaller than 5 
and the same initial data to 
the case $\alpha=1$, $\beta=0.06$, the solution behaves differently. 
The $L^{\infty}$ norm on the left appears to decrease without 
reaching a plateau. The main difference to Fig.~\ref{ben5gauss6em2} is, 
however, in the formation of a strong dispersive shock wave to the 
right of the initial hump. There is still some radiation emanating to 
the left, but it appears that the main part of the initial hump will 
end up in the modulated oscillations propagating towards the right. 
The behavior is similar to what is known from the Kawahara equations, 
a fifth order KdV-type equation, see for instance \cite{DGK} and 
references therein. The dispersive term has the Fourier symbol 
$ik(\alpha |k|-\beta k^{2})$. For $k>0$, plane waves with wave number 
$k<\alpha/\beta$ travel to the left, whereas those with 
$k>\alpha/\beta$ have positive phase velocity. 

\begin{figure}[!htb]
\includegraphics[width=0.49\hsize]{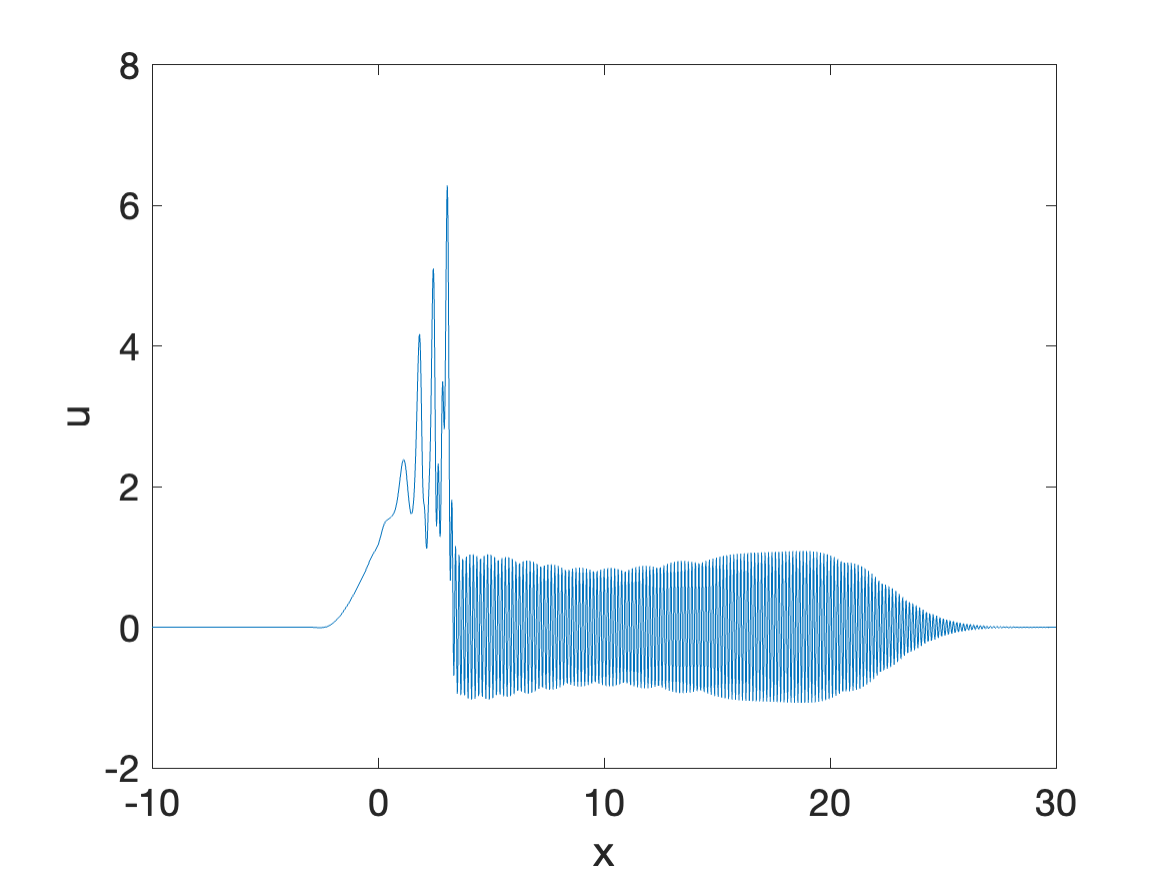}
\includegraphics[width=0.49\hsize]{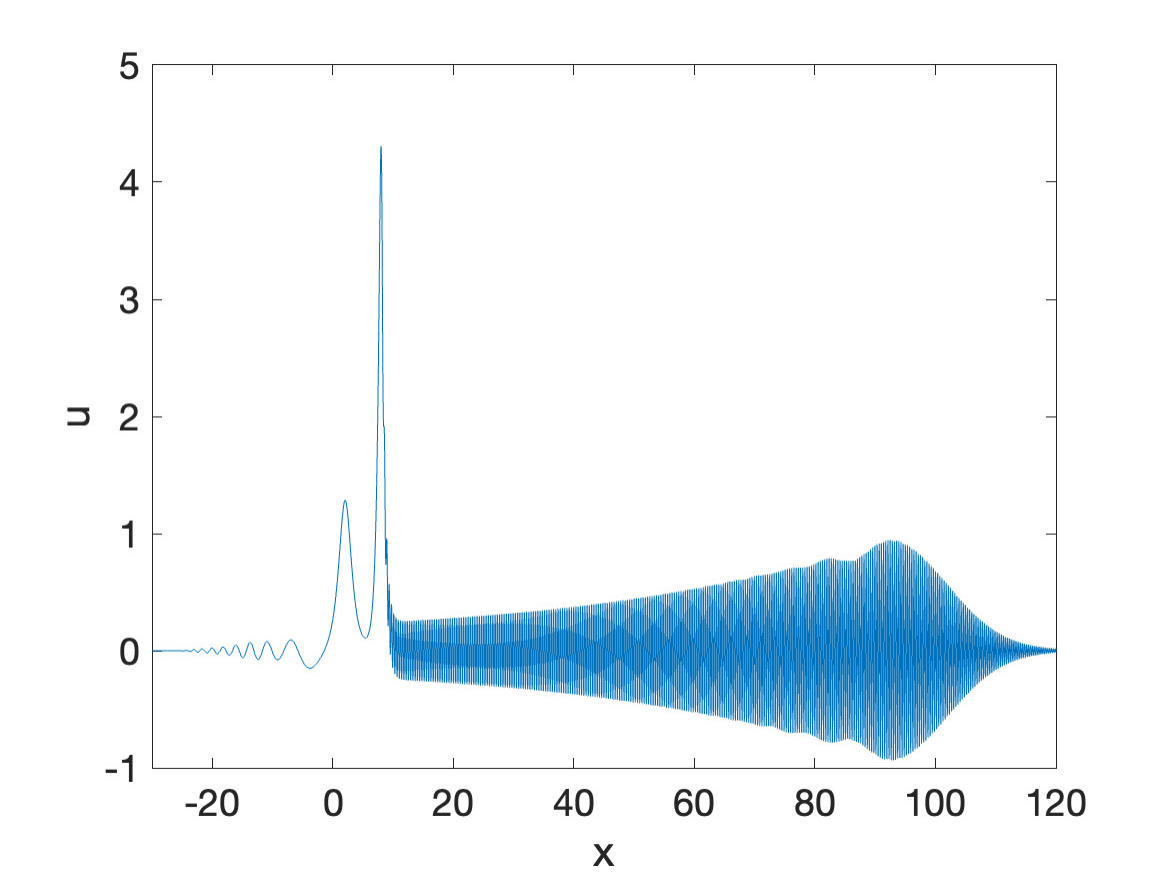}
\caption{Solution to the KdV-ILW equation (\ref{ILWTS})
with $\alpha=1$ for initial data 
$u_{0}(x)=5\exp(-x^{2})$, on the left for $\delta=0.1$ and $t=1$, on 
the right for $\delta=0.9$ and $t=5$.  }
\label{KdVILW5gauss}
\end{figure}

The situation is similar for the KdV-ILW equation (\ref{ILWTS}). The 
solitary waves appear to be stable. Since the figures are very 
similar, we concentrate on the case shown in 
Fig.~\ref{ben5gauss6em2}. When the dispersion and thus the phase 
velocity has different signs for different wave numbers, one gets a 
similar behavior as for the Benjamin equation. We use the same 
numerical parameters and the same initial data and get similar 
solutions in Fig.~\ref{KdVILW5gauss}. Once more there are oscillations 
traveling to the left and rapid oscillations to the right. This 
explains to an extent why we could not find solitary waves in this 
regime. 

\section{Outlook}

\begin{merci}
 F.L. was partially supported by CNPq and FAPERJ/Brazil. C.K. and 
 J.-C.S. were partially supported by the ANR project ISAAC - ANR-23-CE40-0015-01. C.K. 
 thanks for support by  the ANR-17-EURE-0002 EIPHI and by the 
European Union Horizon 2020 research and innovation program under the 
Marie Sklodowska-Curie RISE 2017 grant agreement no. 778010 IPaDEGAN. 
D.P. was supported by a Trond Mohn Foundation grant.

\end{merci}

\end{document}